\documentclass[12pt]{jfm}
\usepackage[usenames,dvipsnames,svgnames,table]{xcolor}  
\usepackage{amsfonts}
\usepackage{fancyhdr}
\usepackage[parfill]{parskip} 
\usepackage{graphicx}
\usepackage{amssymb}
\usepackage{amsmath}
\usepackage{epstopdf}
\usepackage{natbib}
\usepackage{bm,footnpag}
\usepackage{layout}
\usepackage{subcaption}
\usepackage{caption}
\usepackage{cleveref}
\usepackage{float}
\usepackage{xcolor,rotating,picture,lipsum}
\usepackage{wrapfig}
\usepackage{lscape}
\usepackage{rotating}
\usepackage{epstopdf}
\usepackage{multirow}
\usepackage{ifthen}
\usepackage{pgf}
\usepackage{pgffor}
\usepackage[inline]{asymptote} 
\usepackage{parskip}
\usepackage{caption}
\usepackage{subcaption}
\usepackage{pgfplots}
\usepackage{tkz-euclide}
\usepackage{subfig}
\usepackage{placeins}

\ifCUPmtlplainloaded \else
  \checkfont{eurm10}
  \iffontfound
   \IfFileExists{upmath.sty}
      {\typeout{^^JFound AMS Euler Roman fonts on the system,
                   using the 'upmath' package.^^J}%
       \usepackage{upmath}}
      {\typeout{^^JFound AMS Euler Roman fonts on the system, but you
                   dont seem to have the}%
       \typeout{'upmath' package installed. JFM.cls can take advantage
                 of these fonts,^^Jif you use 'upmath' package.^^J}%
      }
  \else
  \fi
\fi

\ifCUPmtlplainloaded \else
  \checkfont{msam10}
  \iffontfound
    \IfFileExists{amssymb.sty}
      {\typeout{^^JFound AMS Symbol fonts on the system, using the
                'amssymb' package.^^J}%
       \usepackage{amssymb}%

      }{}
  \fi
\fi

\ifCUPmtlplainloaded \else
  \IfFileExists{amsbsy.sty}
    {\typeout{^^JFound the 'amsbsy' package on the system, using it.^^J}%
     \usepackage{amsbsy}}
    {}
\fi

\newsavebox{\astrutbox}
\sbox{\astrutbox}{\rule[-5pt]{0pt}{20pt}}

\title[The Navier--Stokes Existence and Smoothness Poser in $\mathbb{R}^n$]{The Navier--Stokes Existence and Smoothness Poser in $\mathbb{R}^n$}
\author[R. K. Michael Thambynayagam]%
{R. K. Michael Thambynayagam%
  \thanks{arXiv:1509.08766 [math.AP], Email: michael.thambynayagam@gmail.com}}
\affiliation{Schlumberger,  Houston, Texas, USA}
  \numberwithin{equation}{section}
\begin{document}
\maketitle
\begin{abstract}
In this paper we describe a method to derive solutions of the incompressible Navier-Stokes system of equations for non-stationary initial  value problems in $\mathbb{R}^n$. We show that for a given smooth solenoidal initial velocity vector field there exist smooth spatially periodic solutions of pressure and velocity in $\mathbb{R}^n$. An illustrative example in $\mathbb{R}^3$ provides important insights into the ostensible phenomenon of the \protect{\itshape{blowup time}}.
\end{abstract}
\protect{\bfseries{2000 Math Subject Classification:}} 76D05, 76D07, 76N10

\section{Introduction}
The object of this paper is to draw attention to a class of problems in which physically reasonable solutions of the non-stationary Navier-Stokes equations may be obtained for incompressible viscous fluids filling all of $\mathbb{R}^n$. 
We recast the Navier-Stokes equation for velocity in terms of three distinct terms associated, respectively, with the linear viscous force, the externally applied force and the inertial force.\\\\
We begin by showing that the two dimensional solution \cite{Tay, Tay2, lad} and a three parameter stationary problem in $\mathbb{R}^3$ introduced by \cite{arn} belong to the case where the components of the term associated with the inertial force, denoted $\mathcal{U}_i\left( x,t \right)$, are zero. We then present a solution of the complete Navier-Stokes equations of motion in $\mathbb{R}^n$ for the case where the components of $\mathcal{U}_i\left( x,t \right)$ are nonzero. An illustrative example in $\mathbb{R}^3$ is included.
\vspace{-0.15in}
\section{The fundamental problem}
The Navier-Stokes equations for viscous incompressible fluids are given by \cite{pdnr}:
%
\begin{eqnarray}
\label{1}
\frac{{\partial v_i }}
{{\partial t}} + g_i= \kappa \Delta v_i -\frac{1}
{\rho }\frac{{\partial p}}
{{\partial x_i }} + f_i,\quad1 \leqslant i \leqslant n,\,\, {x \in \mathbb{R}^n},\quad t \geqslant 0\quad 
\end{eqnarray}
\vspace{-0.20in}
%
\begin{eqnarray}
\label{2}
\protect{\bf{div}}\,\bm{v} = \sum\limits_{i = 1}^n {\frac{{\partial v_i }}
{{\partial x_i }}}=0, \qquad
x \in \mathbb{R}^n,\,\,t \geqslant 0
\end{eqnarray}
where ${v_i} \equiv {v_i}\left( {x,t} \right)$ are the components of velocity vector field at time $t \geqslant 0$, 
%
\begin{eqnarray}
\label{3}
g_i  = \sum_{j = 1}^n {v_j \frac{{\partial v_i }}
{{\partial x_j }}} 
\end{eqnarray}
are the components of the nonlinear inertial force, $p$ is the pressure field, $f_i\left( {x,t} \right)$ are the components of an externally applied force, $\rho$ is the constant density of the fluid,  $\kappa$ is the positive coefficient of kinematical viscosity and $\Delta  = \sum_{i = 1}^n {\frac{{\partial ^2 }}{{\partial x_i^{2} }}}$ is the Laplacian in the space variables.\\\\ The initial conditions are
%
\begin{eqnarray}
\label{4}
v_i\left( {x,0} \right) = v_i^{(0)} 
 \left( x \right),\qquad
x \in \mathbb{R}^n 
\end{eqnarray}
$v_i^{(0)} \left( x \right)$ is a given solenoidal  vector field on $\mathbb{R}^n$. Henceforth, the superscript $(0)$ is used to denote the value of a function at time zero. It is important to note that prescribing pressure at the initial time independent of velocity would render the problem ill-posed.\\\\
Conservation law implies that the energy dissipation of a viscous fluid is bounded by the initial kinetic energy which is finite.  Therefore, if the externally applied force does no net work on the fluid,  the solution must satisfy
%
\begin{eqnarray}
\label{5}
\int\limits_{\mathbb{R}^n } {\left| \bm{v} \right|} ^2 dx < \mathcal{C},\qquad t \geqslant 0
\end{eqnarray}
Reference may be made to \cite{bat} for a formal derivation of the Navier-Stokes system of eqauations. Assuming that $v_i\left( {x,t} \right)$ and $f_i\left( {x,t} \right)$ are smooth and separable in $x$ and $t$, we seek a solution of the Navier-Stokes system of equations $\left(\ref{1}\right)-\left(\ref{4}\right)$ that are spatially periodic in $\mathbb{R}^n$.
%
%
%
\vspace{-0.10in}
\section{Recasting the Navier-Stokes equation: $\mathbb{R}^n=\left\{ { - \infty  < x_i  < \infty ; \,\, i = 1,2,...,n} \right\}$}
Assuming that the divergence and the linear operator can be commuted, the pressure field can be formally obtained by taking the divergence of  $\left(\ref{1}\right)$ as a solution of the Poisson equation, which is
%
\begin{eqnarray}
\label{7}
\Delta p = \rho \sum\limits_{i = 1}^n {\frac{{\partial \left( {f_i  - g_i } \right)}}
{{\partial x_i }}} 
\end{eqnarray}
Equation $\left(\ref{7}\right)$ is called the simplified pressure Poisson equation (PPE). The use of PPE in solving the Navier-Stokes equation is discussed in a paper by \cite{gre}. It is important to note that  while  $\left(\ref{1}\right)$ and $\left(\ref{2}\right)$ lead to the pressure Poisson equation $\left(\ref{7}\right)$, the reverse; that is,  $\left(\ref{1}\right)$ and $\left(\ref{7}\right)$, do not always lead to $\left(\ref{2}\right)$.\\\\ The general solution of the Poisson equation $\left(\ref{7}\right)$ is
%
%
\begin{eqnarray}
\label{8}
p\left( {x,t} \right) &=&- \frac{\rho }
{{2\pi }}\int\limits_{\mathbb{R}^2 } {{\rm P}\left( {y ,t} \right)\ln \left( {\frac{1}
{{\sqrt {\mathcal{P}_n\left( {x,y} \right)} }}} \right)\prod\limits_{j = 1}^2 {dy_j } },\,\, n=2, \nonumber\\
p\left( {x,t} \right) &=&- \frac{{\rho \Gamma \left( {\frac{n}
{2}} \right)}}
{{2\left( {n - 2} \right)\pi ^{\frac{n}
{2}} }}\int\limits_{\mathbb{R}^n } {\frac{{{\rm P}\left( {y,t} \right)}}
{{\left\{ {\mathcal{P}_n\left( {x,y} \right)} \right\}^{\frac{{n - 2}}
{2}} }}} \prod\limits_{j = 1}^n {dy_j },\,\, n\geqslant 3\,\,\, 
\end{eqnarray}
where $\Gamma \left( z \right) = \int_0^\infty  {e^{ - u} u^{z - 1} du}$\quad$\left[ {\Re z > 0} \right]$, is the Gamma function,
%
\begin{eqnarray}
\label{9}
{\rm P}\left( {x,t} \right) = \sum\limits_{j = 1}^n {\frac{{\partial \left( {f_j  - g_j } \right)}}
{{\partial x_j }}} 
\end{eqnarray}
and
%
%
\begin{eqnarray}
\label{10}
\mathcal{P}_n\left( {x,y} \right) = \sum\limits_{j = 1}^n {\left( {x_j  - y_j } \right)^2 } 
\end{eqnarray}
Differentiating $\left(\ref{8}\right)$ with respect to $x_i$ we get 
%
%
\begin{eqnarray}
\label{11}
\frac{{\partial p }}
{{\partial x_i }} = \frac{{\rho \Gamma \left( {\frac{n}
{2}} \right)}}
{{2\pi ^{\frac{n}
{2}} }}\int\limits_{\mathbb{R}^n } {\frac{{\left( {x_i  - y_i } \right){\rm P}\left( {y,t} \right)}}
{{\left\{ {\mathcal{P}_n\left( {x,y } \right)} \right\}^{\frac{n}
{2}} }}} \prod\limits_{j = 1}^n {dy_j },\,\, n\geqslant 2\qquad  
\end{eqnarray}
Substituting for $\frac{{\partial p}}{{\partial x_i }}$ in  $\left(\ref{1}\right)$ we get
%
\begin{eqnarray}
\label{12}
\frac{{\partial v_i }}
{{\partial t}} &=&\kappa \Delta v_i  -\frac{{\Gamma \left( {\frac{n}
{2}} \right)}}
{{2\pi ^{\frac{n}
{2}} }}\int\limits_{\mathbb{R}^n } {\frac{{\left( {x_i  - y_i } \right){\rm P}\left( {y,t} \right)}}
{{\left\{ {\mathcal{P}_n\left( {x,y} \right)} \right\}^{\frac{n}
{2}} }}} \prod\limits_{j = 1}^n {dy_j }+f_i - g_i,\quad n\geqslant 2,\,\, {x \in \mathbb{R}^n},\,\, t \geqslant 0\qquad  
\end{eqnarray}
By using $\left(\ref{9}\right)$, we recast the Navier-Stokes equation $\left(\ref{12}\right)$ as:
%
\begin{eqnarray}
\label{13}
\frac{{\partial v_i }}
{{\partial t}} &=&\kappa \Delta v_i +\mathcal{F}_i\left( {x ,t} \right)  - \mathcal{U}_i\left( x,t \right),\quad {x \in \mathbb{R}^n},\quad t \geqslant 0\quad
\end{eqnarray}
where 
\begin{eqnarray}
\label{14}
\mathcal{F}_i\left( {x ,t} \right) = f_i  - \frac{{\Gamma \left( {\frac{n}
{2}} \right)}}
{{2\pi ^{\frac{n}
{2}} }}\int\limits_{\mathbb{R}^n } {\frac{{\left( {x_i  - w_i } \right)\sum\limits_{k = 1}^n {\frac{{\partial f_k \left( {w ,t} \right)}}
{{\partial w_k }}} }}
{{\left\{ {\mathcal{P}_n\left( {x,w} \right)} \right\}^{\frac{n}
{2}} }}} \prod\limits_{j = 1}^n {dw_j },\quad {x \in \mathbb{R}^n},\,\, t \geqslant 0 \qquad
\end{eqnarray}
and
%
\begin{eqnarray}
\label{15}
\mathcal{U}_i\left( x,t \right)  = g_i  - \frac{{\Gamma \left( {\frac{n}
{2}} \right)}}
{{2\pi ^{\frac{n}
{2}} }}\int\limits_{\mathbb{R}^n } {\frac{{\left( {x_i  - y_i } \right)\sum\limits_{k= 1}^n {\frac{{\partial g_k \left( {y,t} \right)}}
{{\partial y_k }}} }}
{{\left\{ {\mathcal{P}_n\left( {x,y} \right)} \right\}^{\frac{n}
{2}} }}} \prod\limits_{j = 1}^n {dy_j },\quad {x \in \mathbb{R}^n},\,\, t \geqslant 0 \qquad
\end{eqnarray}
We make the following assertions:\\\\
$\bm{\left(i\right)}$ The three terms on the right hand side of  $\left(\ref{13}\right)$, $\kappa \Delta v_i $, $\mathcal{F}_i\left( {x ,t} \right)$ and $\mathcal{U}_i\left( x,t \right)$ are associated, respectively, with the linear viscous force, the externally applied force and the nonlinear inertial force.\\\\ 
$\bm{\left(ii\right)}$ If $\mathcal{U}_i^{(0)} \left( x\right) \equiv 0$; that is
%
\begin{eqnarray}
\label{16}
g_i^{(0)}\left( {x} \right) = \frac{{\Gamma \left( {\frac{n}
{2}} \right)}}
{{2\pi ^{\frac{n}
{2}} }}\int\limits_{\mathbb{R}^n } {\frac{{\left( {x_i  - y_i } \right)\sum\limits_{k = 1}^n {\frac{{\partial g_k^{(0)} }}
{{\partial y_k }}} }}
{{\left\{ {\mathcal{P}_n\left( {x,y} \right)} \right\}^{\frac{n}
{2}} }}} \prod\limits_{j = 1}^n {dy_j }
\end{eqnarray}
then, the solution of the inhomogeneous diffusion equation \cite{tha2} 
%
\begin{eqnarray}
\label{17}
v_i\left( {x,t} \right) &=&\frac{1}
{{\left( {2\sqrt {\pi \kappa t} } \right)^n }}\int\limits_{\mathbb{R}^n } {v_i^{(0)} \left( {y} \right)e^{ - \sum\limits_{k = 1}^n {\frac{{\left( {x_k - y_k } \right)^2 }}
{{4\kappa t}}} } \prod\limits_{j = 1}^n {dy_j } }  + \nonumber\\
&+&\frac{1}{\left( {2\sqrt {\pi \kappa } } \right)^n }\int\limits_{\mathbb{R}^n } {\int\limits_0^t {\frac{{\mathcal{F}_i \left( {y ,\tau } \right)e^{ - \sum\limits_{k = 1}^n {\frac{{\left( {x_k  - y_k } \right)^2 }}
{{4\kappa \left( {t - \tau } \right)}}} } }}
{{\left( {t - \tau } \right)^{\frac{n}{2}}}}d\tau\prod\limits_{j = 1}^n {dy_j } } }
\end{eqnarray}
is a solution of the Navier-Stokes equation $\left(\ref{13}\right)$.\\\\
We now verify assertion $\bm{\left(ii\right)}$ by first deriving two known solutions in $\mathbb{R}^2$ and $\mathbb{R}^3$. \\
\begin{center}
{\protect{\bfseries{The Taylor solution in $\mathbb{R}^2=\left\{ { - \infty  < x_i  < \infty ;\,  i = 1,2} \right\}$\\}}}
\end{center}
\vspace{0.05 in}
$v_1^{(0)}  = \sin \left( {\pi x_1 } \right)\cos \left( {\pi x_2 } \right)$,  $v_2^{(0)}  =  - \cos \left( {\pi x_1 } \right)\sin \left( {\pi x_2 } \right)$ and $f_1=f_2=0$.\\\\ Substituting for $v_1^{(0)}$ and $v_2^{(0)}$ in $\left(\ref{3}\right)$ and $\left(\ref{16}\right)$ we get
%
\begin{eqnarray}
\label{18}
g_i^{(0)}=\frac{1}
{{2\pi }}\int\limits_{\mathbb{R}^2 } {\frac{{\left( {x_i  - y_i } \right)\sum\limits_{k = 1}^2 {\frac{{\partial g_k^{(0)} \left( {y,t} \right)}}
{{\partial y_k}}} }}
{{\mathcal{P}_2 \left( {x,y} \right)}}} \prod\limits_{j = 1}^n {dy_j }  = \frac{\pi }
{2}\sin \left( {2\pi x_i } \right)
\end{eqnarray}
resulting in ${\mathcal{U}_i}\left( {x,t} \right) \equiv 0,\,\, t \geqslant 0$\protect{\footnotemark}. \protect{\footnotetext{In evaluating the integrals in  $\left(\ref{18}\right)$, we have used the following identities:\\ \protect{$\int\limits_{ - \infty }^\infty\!\!  {\sin \left( {\pi u} \right)\!e^{ - \frac{{\left( {x - u} \right)^2 }}
{{4\tau}}} } du = 2\sqrt {\pi \tau} e^{ - \pi ^2 \tau} \sin \left( {\pi x} \right)$ and $\!\int\limits_{ - \infty }^\infty \!\! {\cos \left( {\pi u} \right)\!e^{ - \frac{{\left( {x - u} \right)^2 }}
{{4\tau}}} } du = 2\sqrt {\pi \tau} e^{ - \pi ^2 \tau} \cos \left( {\pi x} \right)$.}}}We obtain $v\left( {x,t} \right)$, $p\left( {x,t} \right)$ from $\left(\ref{17}\right)$ and $\left(\ref{8}\right)$:
%
\begin{eqnarray}
\label{19}
v_1  = \sin \left( {\pi x_1 } \right)\cos \left( {\pi x_2 } \right)e^{ - 2\pi ^2 \kappa t}
\end{eqnarray}
%
\begin{eqnarray}
\label{20}
v_2  =  - \cos \left( {\pi x_1 } \right)\sin \left( {\pi x_2 } \right)e^{ - 2\pi ^2 \kappa t}
\end{eqnarray}
and
\begin{eqnarray}
\label{21}
p =  - \frac{{\rho e^{ - 4\pi ^2 \kappa t} }}{4}\left[ {\cos \left( {2\pi x_1 } \right) + \cos \left( {2\pi x_2 } \right)} \right]
\end{eqnarray}
which is the two dimensional \cite{Tay} solution on the exponential decay of vortices in a viscous fluid.\\
\begin{center}
{\protect{\bfseries{Arnold-Beltrami-Childress (ABC) flows in $\mathbb{R}^3=\left\{ { - \infty  < x_i  < \infty ; \,\, i = 1,2,3} \right\}$\\}}}
\end{center}
\vspace{0.1 in}
$v_1^{(0)}  = a\sin \pi x_3  - c\cos \pi x_2$, $v_2^{(0)}  = b\sin \pi x_1 - a\cos \pi x_3$, $v_3^{(0)}  = c\sin \pi x_2  - b\cos \pi x_1$
 and $f_1=f_2=f_3=0$. $a$, $b$ and $c$ are real constants.\\\\ Substituting for $v_1^{(0)}$, $v_2^{(0)}$ and $v_3^{(0)}$ in $\left(\ref{3}\right)$ and $\left(\ref{16}\right)$ we obtain
%
\begin{eqnarray}
\label{22}
\!\!\!\!\!\!\!\!g_1^{(0)}  &=& \frac{1}
{{4\pi }}\int\limits_{\mathbb{R}^3 } {\frac{{\left( {x_1  - y_1 } \right)\sum\limits_{k = 1}^3 {\frac{{\partial g_k^{(0)} \left( {y ,t} \right)}}
{{\partial y_k }}} }}
{{\left\{ {\mathcal{P}_3\left( {x,y} \right)} \right\}^{\frac{3}
{2}} }}} \prod\limits_{j = 1}^3 {dy_j }
\nonumber\\
&  = & 
\pi \left\{ {bc\sin \left(\pi x_1\right) \sin \left(\pi x_2\right)- ab\cos \left(\pi x_1\right) \cos\left( \pi x_3\right) } \right\}\qquad
\end{eqnarray}
%
\begin{eqnarray}
\label{23}
\!\!\!\!\!\!\!\!g_2^{(0)}  &  = & \frac{1}
{{4\pi }}\int\limits_{\mathbb{R}^3 } {\frac{{\left( {x_2  - y_2 } \right)\sum\limits_{k = 1}^3 {\frac{{\partial g_k^{(0)} \left( {y ,t} \right)}}
{{\partial y_k }}} }}
{{\left\{ {\mathcal{P}_3\left( {x,y} \right)} \right\}^{\frac{3}
{2}} }}} \prod\limits_{j = 1}^3 {dy_j } 
\nonumber\\
&  = &\pi \left\{ {ac\sin \left(\pi x_2\right) \sin \left(\pi x_3\right)-bc\cos \left(\pi x_1\right) \cos \left(\pi x_2\right) } \right\}\qquad
\end{eqnarray}
%
\begin{eqnarray}
\label{24}
\!\!\!\!\!\!\!\!g_3^{(0)} &  = &\frac{1}
{{4\pi }}\int\limits_{\mathbb{R}^3 } {\frac{{\left( {x_3  - y_3 } \right)\sum\limits_{k = 1}^3 {\frac{{\partial g_k^{(0)} \left( {y,t} \right)}}
{{\partial y_k}}} }}
{{\left\{ {\mathcal{P}_3\left( {x,y } \right)} \right\}^{\frac{3}
{2}} }}} \prod\limits_{j = 1}^3 {dy_j } 
\nonumber\\
&  = & \pi \left\{ {ab\sin\left( \pi x_1\right) \sin \left(\pi x_3\right)-ac\cos\left( \pi x_2\right) \cos \left(\pi x_3\right)} \right\}\qquad
\end{eqnarray}
which results in ${\mathcal{U}_i}\left( {x,t} \right) \equiv 0,\,\, t \geqslant 0$. We obtain $v\left( {x,t} \right)$, $p\left( {x,t} \right)$ from $\left(\ref{17}\right)$ and $\left(\ref{8}\right)$:
%
\begin{eqnarray}
\label{25}
v_1 = \left\{ {a\sin \left( {\pi x_3 } \right) - c\cos \left( {\pi x_2 } \right)} \right\}e^{ - \pi ^2 \kappa t} 
\end{eqnarray}
%
\begin{eqnarray}
\label{26}
v_2  = \left\{b\sin \left( {\pi x_1 } \right)- {a\cos \left( {\pi x_3 } \right)} \right\}e^{ - \pi ^2 \kappa t} 
\end{eqnarray}
%
\begin{eqnarray}
\label{27}
v_3  = \left\{ {c\sin \left( {\pi x_2 } \right) - b\cos \left( {\pi x_1 } \right)} \right\}e^{ - \pi ^2 \kappa t} 
\end{eqnarray}
and
\begin{eqnarray}
\label{28}
p \!= \! - \rho e^{ - 2\pi ^2 \kappa t}\! \left[ {bc\cos \left( {\pi x_1 } \right)\sin \left( {\pi x_2 } \right) + } {ab\cos \left( {\pi x_3 } \right)\sin \left( {\pi x_1 } \right) + ac\cos \left( {\pi x_2 } \right)\sin \left( {\pi x_3 } \right)} \right]\quad\quad
\end{eqnarray}
The three parameter ($a$, $b$ and $c$) family of periodic flows provide a simple non-stationary solution of the Navier-Stokes equation which is the exponential decay of helical streamlines. The initial  solution ($t=0$) of which, however, is also the stationary solution of the inviscid Euler equation independently introduced by \cite{arn} and \cite{chi1,chi2}. The flow has complex characteristics which have been studied exhaustively by \cite{Dom}.
\vspace{-0.20in}
\section{The solution in $\mathbb{R}^n$, $\mathcal{U}_i^{(0)} \left( x\right)  \ne 0$, $i=1,2...,n$}
We set the components of the externally applied force $f_i\left( {x,t} \right)$ to zero and observe that, for a given solenoidal and spatially periodic initial velocity vector field, the nonlinearity \protect{\itshape{instantaneously}} spirals from zero to a plateau through a sequence of linear diffusion processes in accordance with $\left(\ref{13}\right)$. The plateaued ${\mathcal{U}_i}\left( x,t \right)$, substituted into $\left(\ref{13}\right)$, results in a linear inhomogeneous partial differential equation, which is solved analytically. \\\\
\vspace{-0.20in}
\begin{center}
{\protect{\bfseries{The \protect{\itshape{instantaneous}} generation and plateauing of ${\mathcal{U}_i}\left( x,t \right)$}}}
\end{center}
The velocity vector field $v_i^{(0)} \left( {x} \right)$ is a periodic function. Therefore, $\mathcal{U}_i^{(0)} \left( {x} \right)$ is also a periodic function.
We view $\mathcal{U}_i \left( {x,t} \right)$ as a generating function both spatially and in time; that is, for a given positive coefficient of kinematical viscosity, $\kappa$, ${\mathcal{U}_i}\left( x,t \right)$, is manifested through an \protect{\itshape{instantaneous sequence}} $\mathcal{U}_i^{(j)}\left( x, t \right)$, $j=0,1,...$ until it plateaus. Where the sequence is denoted by the superscript $(j)$.  When $\mathcal{U}_i^{(k)}\left( x, t \right) = \mathcal{U}_i^{(k - 1)}\left( x, t \right)$ we conclude that the manifestation of nonlinearity has ceased and set ${\mathcal{U}_i}\left( x,t \right) = \mathcal{U}_i^{(k - 1)}\left( x, t \right)$ in $\left(\ref{13}\right)$. The solution of $\left(\ref{13}\right)$ at sequence $(k-1)$ satisfies the Navier-Stokes system of equations $\left(\ref{1}\right)-\left(\ref{4}\right)$.\\\\
The concept is best explained by formally deriving the explicit formula for the generating function ${\mathcal{U}_i}\left( x,t \right)$. For the sake of simplicity, we choose a velocity vector field that exhibits spatial symmetry.\\\\
{\protect{\bfseries{Sequence: 1}}}\\\\
We begin the sequencing with no nonlinearity by setting ${\mathcal{U}_i}\left( x,t \right)=0$ in $\left(\ref{13}\right)$. The first sequence velocity vector field is given by
\begin{eqnarray}
\label{29}
v_i^{(1)}\left( x,t \right) = \frac{1}{{{{\left( {2\sqrt {\pi \kappa t} } \right)}^n}}}\int\limits_{{\mathbb{R}^n}} {v_i^{(0)}\left( y \right){e^{ - \sum\limits_{k = 1}^n {\frac{{{{\left( {{x_k} - {y_k}} \right)}^2}}}{{4\kappa t}}} }}\prod\limits_{j = 1}^n {d{y_j}} }
\end{eqnarray}
Making use of  the integral identities 
%
\begin{eqnarray}
\label{30}
 \int\limits_{ - \infty }^\infty  {\sin \left( {\alpha u} \right){e^{ - \frac{{{{\left( {x - u} \right)}^2}}}{{4\tau }}}}} du = 2\sqrt {\pi \tau } {e^{ - {\alpha ^2}\tau }}\sin \left( {\alpha x} \right) 
\end{eqnarray}
%
\begin{eqnarray}
\label{31}
 \int\limits_{ - \infty }^\infty  {\cos \left( {\alpha u} \right){e^{ - \frac{{{{\left( {x - u} \right)}^2}}}{{4\tau }}}}} du = 2\sqrt {\pi \tau } {e^{ - {\alpha ^2}\tau }}\cos \left( {\alpha x} \right)
\end{eqnarray}
and
%
\begin{eqnarray}
\label{32}
\int\limits_{ - \infty }^\infty  {{e^{ - \frac{{{{\left( {x - u} \right)}^2}}}{{4\tau }}}}} du = 2\sqrt {\pi \tau }
\end{eqnarray}
the right hand side of $\left(\ref{29}\right)$ may be written as
\begin{eqnarray}
\label{33}
\frac{1}{{{{\left( {2\sqrt {\pi \kappa t} } \right)}^n}}}\int\limits_{{\mathbb{R}^n}} {v_i^{(0)}\left( y \right){e^{ - \sum\limits_{k = 1}^n {\frac{{{{\left( {{x_k} - {y_k}} \right)}^2}}}{{4\kappa t}}} }}\prod\limits_{j = 1}^n {d{y_j}} }  =  v_i^{(0)}{e^{ - {\xi_1}{\pi ^2}\kappa t}}
\end{eqnarray}
We have
%
\begin{eqnarray}
\label{34}
v_i^{(1)}\left( x,t \right) = v_i^{(0)}\left( x \right){\mathcal{T} _1}\left( t \right)
\end{eqnarray}
Where
%
\begin{eqnarray}
\label{35}
{\mathcal{T} _1}\left( t \right) = {e^{ - {\xi_1}{\pi ^2}\kappa t}}
\end{eqnarray}
and $\xi_1$ is a real positive constant resulting from performing the integrals over the initial velocity vector field $v_i^{(0)}\left( x \right)$ comprising circular functions in $\mathbb{R}^n$.\\\\
We obtain $\mathcal{U}_i^{(0)}\left( x \right)$ from $\left(\ref{15}\right)$:
%
\begin{eqnarray}
\label{36}
\mathcal{U}_i^{(0)}\left( x \right) = g_i^{(0)}\left( x \right) - \hbar _i^{(0)}\left( x \right)
\end{eqnarray}
where $g_i^{(0)}\left( x\right) $ and $\hbar _i^{(0)}\left( x \right) $ are given by
%
\begin{eqnarray}
\label{37}
g_i^{(0)}\left( x \right) = \sum\limits_{j = 1}^n {v_j^{(0)} } \frac{{\partial v_i^{(0)} }}{{\partial {x_j}}}
\end{eqnarray}
and
%
\begin{eqnarray}
\label{38}
\hbar _i^{(0)}\left( x \right) = \frac{{\Gamma \left( {\frac{n}{2}} \right)}}{{2{\pi ^{\frac{n}{2}}}}}\int\limits_{{\mathbb{R}^n}} {\frac{{\left( {{x_i} - {y_i}} \right)\sum\limits_{k = 1}^n {\frac{{\partial g_k^{(0)}\left( {y} \right)}}{{\partial {y_k}}}} }}{{{{\left\{ {{\mathcal{P}_n}\left( {x,y} \right)} \right\}}^{\frac{n}{2}}}}}} \prod\limits_{j = 1}^n {d{y_j}} 
\end{eqnarray}
{\protect{\bfseries{Sequence: 2}}}
%
\begin{eqnarray}
\label{39}
\mathcal{U}_i^{(1)}\left( x,t \right) = g_i^{(1)}\left( x,t \right)  - \hbar _i^{(1)}\left( x,t \right) 
\end{eqnarray}
where $g_i^{(1)}\left( x,t \right) $ and $\hbar _i^{(1)}\left( x,t \right) $ are given by
%
\begin{eqnarray}
\label{40}
g_i^{(1)}\left( {x,t} \right) = \sum\limits_{j = 1}^n {v_j^{(1)}} \frac{{\partial v_i^{(1)}}}{{\partial {x_j}}} = g_i^{(0)}\left( x \right)\mathcal{T} _1^2\left( t \right) = g_i^{(0)}\left( x \right) {e^{ - 2{\xi _1}{\pi ^2}\kappa t}}
\end{eqnarray}
and
%
\begin{eqnarray}
\label{41}
\hbar _i^{(1)}\left( x,t \right)  = \frac{{\Gamma \left( {\frac{n}{2}} \right)}}{{2{\pi ^{\frac{n}{2}}}}}\int\limits_{{\mathbb{R}^n}} {\frac{{\left( {{x_i} - {y_i}} \right)\sum\limits_{k = 1}^n {\frac{{\partial g_k^{(1)}\left( {y,t} \right)}}{{\partial {y_k}}}} }}{{{{\left\{ {{\mathcal{P}_n}\left( {x,y} \right)} \right\}}^{\frac{n}{2}}}}}} \prod\limits_{j = 1}^n {d{y_j}}  = \hbar _i^{(0)}\left( x \right) {e^{ - 2{\xi_1}{\pi ^2}\kappa t}}
\end{eqnarray}
respectively. Substituting for $g_i^{(1)}\left( x,t \right) $ and $\hbar _i^{(1)}\left( x,t \right) $ in $\left(\ref{39}\right)$, we obtain
%
\begin{eqnarray}
\label{42}
\mathcal{U}_i^{(1)}\left( x,t \right)  = \mathcal{U}_i^{(0)}\left( x \right){e^{ - 2{\xi_1}{\pi ^2}\kappa t}}
\end{eqnarray}
Setting ${\mathcal{U}_i}\left( x,t \right) = \mathcal{U}_i^{(1)}\left( x,t \right)  = \mathcal{U}_i^{(0)}\left( x \right){e^{ -2{\xi_1}{\pi ^2}\kappa t}}$ in $\left(\ref{13}\right)$ and solving the inhomogeneous diffusion equation, we obtain:
%
\begin{eqnarray}
\label{43}
v_i^{(2)}\left( x,t \right)=v_i^{(1)}\left( x,t \right)-\frac{1}{{{{\left( {2\sqrt {\pi \kappa } } \right)}^n}}}\int\limits_0^t {\int\limits_{{\mathbb{R}^n}} {\frac{{\mathcal{U}_i^{(0)}\left( y \right){e^{ - 2{\xi_1}{\pi ^2}\kappa \tau }}{e^{ - \sum\limits_{k = 1}^n {\frac{{{{\left( {{x_k} - {y_k}} \right)}^2}}}{{4\kappa \left( {t - \tau } \right)}}} }}}}{{{{\left( {t - \tau } \right)}^{\frac{n}{2}}}}}\prod\limits_{j = 1}^n {d{y_j}} } } d\tau \qquad
\end{eqnarray}
Integrating the second term on the right-hand side of $\left(\ref{43}\right)$ we obtain
\begin{eqnarray}
\label{44}
\frac{1}{{{{\left( {2\sqrt {\pi \kappa } } \right)}^n}}}\!\!\int\limits_0^t\! {\int\limits_{{\mathbb{R}^n}} {\frac{{\mathcal{U}_i^{(0)}\left( y \right){e^{ - 2{\xi_1}{\pi ^2}\kappa \tau }}{e^{ - \sum\limits_{k = 1}^n {\frac{{{{\left( {{x_k} - {y_k}} \right)}^2}}}{{4\kappa \left( {t - \tau } \right)}}} }}}}{{{{\left( {t - \tau } \right)}^{\frac{n}{2}}}}}\!\prod\limits_{j = 1}^n {d{y_j}} } } d\tau =\frac{{\mathcal{U}_i^{(0)}\!\!\left( x \right)\!\left( {{e^{ -2{\xi_1}{\pi ^2}\kappa t}}\! - \!{e^{ - {\xi_2}{\pi ^2}\kappa t}}} \right)}}{{{\pi ^2}\kappa\left( {\xi_2 - 2\xi_1} \right) }},\nonumber\\
{\xi _2} \ne 2{\xi _1}\qquad\qquad
\end{eqnarray}
Where $\xi_2$ is a real positive constant resulting from performing the integrals over $\mathcal{U}_i^{(0)}\left( y \right)$ comprising circular functions in $\mathbb{R}^n$.\\\\ Substituting  $\left(\ref{44}\right)$ in $\left(\ref{43}\right)$ we obtain the solution of the second sequence velocity vector field:
%
\begin{eqnarray}
\label{45}
v_i^{(2)}\left( x,t \right)  = v_i^{(1)}\left( x,t \right)  - \mathcal{U}_i^{(0)}\left( x \right){\mathcal{T} _2}\left( t \right)
\end{eqnarray}
where
%
\begin{eqnarray}
\label{46}
{\mathcal{T} _2}\left( t \right) = \frac{{\left( {{e^{ - 2{\xi_1}{\pi ^2}\kappa t}} - {e^{ - {\xi_2}{\pi ^2}\kappa t}}} \right)}}{{{\pi ^2}\kappa \left( {\xi_2 - 2\xi_1} \right)}},\quad {\xi _2} \ne 2{\xi _1}
\end{eqnarray}
{\protect{\bfseries{Sequence: 3}}}
%
\begin{eqnarray}
\label{47}
\mathcal{U}_i^{(2)}\left( x,t \right)  = g_i^{(2)}\left( x,t \right)  - \hbar _i^{(2)}\left( x,t \right) 
\end{eqnarray}
where 
%
\begin{eqnarray}
\label{48}
g_i^{(2)}\left( {x,t} \right) = \sum\limits_{j = 1}^n {v_j^{(2)}} \frac{{\partial v_i^{(2)}}}{{\partial {x_j}}} &=& \alpha _{i0}^{(2)}\left( x \right)\mathcal{T} _1^2\left( t \right) + \alpha _{i1}^{(2)}\left( x \right){\mathcal{T} _1}\left( t \right){\mathcal{T} _2}\left( t \right) + \alpha _{i2}^{(2)}\left( x \right)\mathcal{T} _2^2\left( t \right)\nonumber\\
&=& g_i^{(1)}\left( {x,t} \right) + \alpha _{i1}^{(2)}\left( x \right){\mathcal{T} _1}\left( t \right){\mathcal{T} _2}\left( t \right) + \alpha _{i2}^{(2)}\left( x \right)\mathcal{T} _2^2\left( t \right)
\end{eqnarray}
Sequence 3 has generated, for $g_i^{(2)}\left( x,t \right) $, three spatial functions  augmented by exponentially decaying functions of time.  The coefficients $\alpha _{ik}^{(2)}\left(x\right)$, $k=0,1,2$, which are functions of $v_i^{(0)}\left( x \right)$, $\mathcal{U}_i^{(0)}\left( x \right)$ and their derivatives, are given by
\begin{eqnarray}
\label{49}
\alpha _{i0}^{(2)}\left( x \right) = g_i^{(0)}\left( x \right) 
\end{eqnarray}
%
\begin{eqnarray}
\label{50}
\alpha _{i1}^{(2)}\left(x\right) =   -\sum\limits_{j = 1}^n {\mathcal{U}_j^{(0)}\left( x \right) \frac{{\partial v_i^{(0)}\left( x \right) }}{{\partial {x_j}}}}  - \sum\limits_{j = 1}^n {v_j^{(0)}\left( x \right) \frac{{\partial \mathcal{U}_i^{(0)}\left( x \right) }}{{\partial {x_j}}}} 
\end{eqnarray}
%
\begin{eqnarray}
\label{51}
\alpha _{i2}^{(2)}\left(x\right) = \sum\limits_{j = 1}^n {\mathcal{U}_j^{(0)}\left( x \right) \frac{{\partial \mathcal{U}_i^{(0)}\left( x \right) }}{{\partial {x_j}}}} 
\end{eqnarray}
and
%
\begin{eqnarray}
\label{52}
\hbar _i^{(2)}\left( x,t \right)&=& \frac{{\Gamma \left( {\frac{n}{2}} \right)}}{{2{\pi ^{\frac{n}{2}}}}}\int\limits_{{\mathbb{R}^n}} {\frac{{\left( {{x_i} - {y_i}} \right)\sum\limits_{i = 1}^n {\frac{{\partial g_i^{(2)}\left( {{y_j},t} \right)}}{{\partial {x_i}}}} }}{{{{\left\{ {\mathcal{P}\left( {{y_j},{x_j}} \right)} \right\}}^{\frac{n}{2}}}}}} \prod\limits_{j = 1}^n {d{y_j}}\nonumber\\
& = &\beta _{i0}^{(2)}\left( x \right)\mathcal{T} _1^2\left( t \right) + \beta _{i1}^{(2)}{\mathcal{T} _1}\left( t \right){\mathcal{T} _2}\left( t \right) + \beta _{i2}^{(2)}\mathcal{T} _2^2\left( t \right)\nonumber\\
& = &\hbar _i^{(1)}\left( x,t \right) + \beta _{i1}^{(2)}{\mathcal{T} _1}\left( t \right){\mathcal{T} _2}\left( t \right) + \beta _{i2}^{(2)}\mathcal{T} _2^2\left( t \right)
\end{eqnarray}
It is important that we express the coefficients $\alpha _{ik}^{(2)}\left(x\right)$, $k=0,1,2$ in a form integrable by use of the identities $\left(\ref{30}\right)$--$\left(\ref{32}\right)$. The coefficients $\beta _{ik}^{(2)}\left( x \right)$, $k=0,1,2$ are corollary to $\alpha _{ik}^{(2)}\left( x \right)$ and are given by
%
\begin{eqnarray}
\label{53}
\beta _{i0}^{(2)}\left( x \right)  &=&\hbar _i^{(0)}\left( x \right) \nonumber\\
\beta _{ik}^{(2)}\left( x \right)  &=& \frac{{\Gamma \left( {\frac{n}{2}} \right)}}{{2{\pi ^{\frac{n}{2}}}}}\int\limits_{{\mathbb{R}^n}} {\frac{{\left( {{x_i} - {y_i}} \right)\sum\limits_{i = 1}^n {\frac{{\partial \alpha _{ik}^{(2)}\left( {{y_j}} \right)}}{{\partial {x_i}}}} }}{{{{\left\{ {\mathcal{P}\left( {{y_j},{x_j}} \right)} \right\}}^{\frac{n}{2}}}}}} \prod\limits_{j = 1}^n {d{y_j}}, \qquad k=1,2
\end{eqnarray}
Substituting for $g_i^{(2)}\left( x,t \right) $ and $\hbar _i^{(2)}\left( x,t \right) $ in $\left(\ref{47}\right)$, we have
%
\begin{eqnarray}
\label{54}
\mathcal{U}_i^{(2)}\left( x,t \right)  = \mathcal{U}_i^{(1)}\left( x,t \right) + q_i^{(2)}\left( x,t \right) 
\end{eqnarray}
where 
%
\begin{eqnarray}
\label{55}
q_i^{(2)}\left( {x,t} \right) = \mathcal{G}_{i1}^{(2)}\left( x \right){\mathcal{T} _1}\left( t \right){\mathcal{T} _2}\left( t \right) + \mathcal{G}_{i2}^{(2)}\left( x \right){\mathcal{T} _2^2\left( t \right)}
\end{eqnarray}
%
\begin{eqnarray}
\label{56}
{\mathcal{T} _1}\left( t \right){\mathcal{T} _2}\left( t \right) = \frac{{\left( {{e^{ - 3{\xi_1}{\pi ^2}\kappa t}} - {e^{ - \left( {\xi_1+\xi_2} \right){\pi ^2}\kappa t}}} \right)}}{{\left( {\xi_2 - 2\xi_1} \right){\pi ^2}\kappa }},\quad {\xi _2} \ne 2{\xi _1}
\end{eqnarray}
%
\begin{eqnarray}
\label{57}
{\mathcal{T} _2^2\left( t \right)} = \frac{{\left( {{e^{ - 4{\xi_1}{\pi ^2}\kappa t}} - 2{e^{ - \left( {2\xi_1+\xi_2}\right){\pi ^2}\kappa t}} + {e^{ - 2{\xi_2}{\pi ^2}\kappa t}}} \right)}}{{{{\left( {\xi_2 - 2\xi_1} \right)}^2}{\pi ^4}{\kappa ^2}}},\quad {\xi _2} \ne 2{\xi _1}
\end{eqnarray}
and
%
\begin{eqnarray}
\label{58}
\mathcal{G}_{ik}^{(2)}\left( x \right) = \alpha _{ik}^{(2)}\left( x \right)- \beta _{ik}^{(2)}\left( x \right),\qquad k=1,2 
\end{eqnarray}
Substituting for ${\mathcal{U}_i}\left( x,t \right) = \mathcal{U}_i^{(2)}\left( x,t \right) $ in $\left(\ref{13}\right)$ and solving the inhomogeneous diffusion equation, we obtain:
%
\begin{eqnarray}
\label{59}
v_i^{(3)}\left( x,t \right)  &=& v_i^{(1)}\left( x,t \right)  - \frac{1}{{{{\left( {2\sqrt {\pi \kappa } } \right)}^n}}}\int\limits_0^t {\int\limits_{{\mathbb{R}^n}} {\frac{{\mathcal{U}_i^{(2)}\left( y \right){e^{ - \sum\limits_{k = 1}^n {\frac{{{{\left( {{x_k} - {y_k}} \right)}^2}}}{{4\kappa \left( {t - \tau } \right)}}} }}}}{{{{\left( {t - \tau } \right)}^{\frac{n}{2}}}}}\prod\limits_{j = 1}^n {d{y_j}} } } d\tau  \nonumber\\
&=&  v_i^{(1)}\left( x,t \right) - \frac{1}{{{{\left( {2\sqrt {\pi \kappa } } \right)}^n}}}\int\limits_0^t {\int\limits_{{\mathbb{R}^n}} {\frac{{\mathcal{U}_i^{(1)}\left( y \right){e^{ - \sum\limits_{k = 1}^n {\frac{{{{\left( {{x_k} - {y_k}} \right)}^2}}}{{4\kappa \left( {t - \tau } \right)}}} }}}}{{{{\left( {t - \tau } \right)}^{\frac{n}{2}}}}}\prod\limits_{j = 1}^n {d{y_j}} } } d\tau -\nonumber\\
   &&- \frac{1}{{{{\left( {2\sqrt {\pi \kappa } } \right)}^n}}}\int\limits_0^t {\int\limits_{{\mathbb{R}^n}} {\frac{{q_i^{(2)}\left( {y,t} \right){e^{ - \sum\limits_{k = 1}^n {\frac{{{{\left( {{x_k} - {y_k}} \right)}^2}}}{{4\kappa \left( {t - \tau } \right)}}} }}}}{{{{\left( {t - \tau } \right)}^{\frac{n}{2}}}}}\prod\limits_{j = 1}^n {d{y_j}} } } d\tau \nonumber\\
  &=& v_i^{(2)}\left( x,t \right)  - \mathcal{O}_i^{(2)}\left( {x,t} \right)
\end{eqnarray}
where
%
\begin{eqnarray}
\label{60}
\mathcal{O}_i^{(2)}\left( {x,t} \right)=\frac{1}{{{{\left( {2\sqrt {\pi \kappa } } \right)}^n}}}\int\limits_0^t {\int\limits_{{\mathbb{R}^n}} {\frac{{q_i^{(2)}\left( {y,t} \right){e^{ - \sum\limits_{k = 1}^n {\frac{{{{\left( {{x_k} - {y_k}} \right)}^2}}}{{4\kappa \left( {t - \tau } \right)}}} }}}}{{{{\left( {t - \tau } \right)}^{\frac{n}{2}}}}}\prod\limits_{j = 1}^n {d{y_j}} } } d\tau 
\end{eqnarray}
Expressing $q_i^{(2)}\left( {x,t} \right) $ in a form integrable by use of the identities $\left(\ref{30}\right)$--$\left(\ref{32}\right)$ and performing the integrations in $\left(\ref{60}\right)$ we obtain the solution of the third sequence velocity vector field.\\\\ 
{\protect{\bfseries{Sequence:}}\,\,$l$}\\\\
The prescription for obtaining the $l^{th}$ sequence velocity vector field is as follows:\\\\
$\bm{\left(i\right)}$ Compute 
\begin{center}
$\begin{gathered}
  g_i^{(0)}\left( x \right) = \sum\limits_{j = 1}^n {v_j^{(0)}\left( x \right) } \frac{{\partial v_i^{(0)}\left( x \right) }}{{\partial {x_j}}} \hfill \\
  g_i^{(1)}\left( {x,t} \right) = \sum\limits_{j = 1}^n {v_j^{(1)}\left( x,t \right)} \frac{{\partial v_i^{(1)}\left( x,t \right)}}{{\partial {x_j}}} = g_i^{(0)}\left( x \right) {e^{ - 2{\xi _1}{\pi ^2}\kappa t}} \hfill \\ 
\end{gathered}$
\end{center}
%
\begin{eqnarray}
\label{61}
g_i^{(l)}\left( {x,t} \right) = \sum\limits_{j = 1}^n {v_j^{(l)}} \frac{{\partial v_i^{(l)}}}{{\partial {x_j}}} = g_i^{(l - 1)}\left( {x,t} \right) + \sum\nolimits_{m} {{\rm A}\left( {x,t} \right)},\quad l=2,3,...
\end{eqnarray}
\begin{center}
$\begin{gathered}
  \hbar _i^{(0)}\left( x \right) = \frac{{\Gamma \left( {\frac{n}{2}} \right)}}{{2{\pi ^{\frac{n}{2}}}}}\int\limits_{{\mathbb{R}^n}} {\frac{{\left( {{x_i} - {y_i}} \right)\sum\limits_{k = 1}^n {\frac{{\partial g_k^{(0)}\left( y \right)}}{{\partial {y_k}}}} }}{{{{\left\{ {{\mathcal{P}_n}\left( {x,y} \right)} \right\}}^{\frac{n}{2}}}}}} \prod\limits_{j = 1}^n {d{y_j}}  \hfill \\
  \hbar _i^{(1)}\left( {x,t} \right) = \hbar _i^{(0)}{e^{ - 2{\xi _1}{\pi ^2}\kappa t}} \hfill \\ 
\end{gathered}$
\end{center}
%
\begin{eqnarray}
\label{62}
\hbar _i^{(l)}\left( {x,t} \right)& =& \frac{{\Gamma \left( {\frac{n}{2}} \right)}}{{2{\pi ^{\frac{n}{2}}}}}\int\limits_{{\mathbb{R}^n}} {\frac{{\left( {{x_i} - {y_i}} \right)\sum\limits_{k = 1}^n {\frac{{\partial g_k^{(l)}\left( {{y_j},t} \right)}}{{\partial {x_k}}}} }}{{{{\left\{ {\mathcal{P}\left( {{y_j},{x_j}} \right)} \right\}}^{\frac{n}{2}}}}}} \prod\limits_{j = 1}^n {d{y_j}}\nonumber\\
&  = &\hbar _i^{(l-1)}\left( {x,t} \right) + \sum\nolimits_{m} {{\rm B}\left( {x,t} \right)},\quad l=2,3,...
\end{eqnarray}
where $\sum\nolimits_{m} {{\rm A}\left( {x,t} \right)}$ and $\sum\nolimits_{m}{{\rm B}\left( {x,t} \right)}$ are additional $m$ terms generated at sequence $l$.
%
\begin{eqnarray}
\label{63}
\mathcal{U}_i^{(0)}\left( x \right) &=& g_i^{(0)}\left( x \right) - \hbar _i^{(0)}\left( x \right)\nonumber\\
\mathcal{U}_i^{(1)}\left( x,t \right) &=& g_i^{(1)}\left( x,t \right)  - \hbar _i^{(1)}\left( x,t \right) = \mathcal{U}_i^{(0)}\left( x \right){e^{ - 2{\xi_1}{\pi ^2}\kappa t}}\nonumber\\
\mathcal{U}_i^{(l)}\left( {x,t} \right) &=& g_i^{(l)}\left( {x,t} \right) - \hbar _i^{(l)}\left( {x,t} \right) = \mathcal{U}_i^{(l - 1)}\left( {x,t} \right) + q_i^{(l)}\left( {x,t} \right),\quad l=2,3,...
\end{eqnarray}
\begin{center}
$\begin{gathered}
\bm{\left(ii\right)}\,\,\protect{\textrm{By definition}}\quad \mathcal{O}_i^{(0)}\left( {x} \right) = 0\qquad\qquad\qquad\qquad\qquad\qquad\qquad\qquad\qquad\qquad\qquad\qquad\qquad \hfill \\
\qquad \qquad \qquad \qquad \quad \mathcal{O}_i^{(1)}\left( {x,t} \right) = \mathcal{U}_i^{(0)}\left( x \right){\mathcal{T} _2}\left( t \right) \hfill \\ 
\end{gathered}$
\end{center}
Express the coefficients $q_i^{(l)}\left( {x,t} \right)$ in a form integrable by use of the identities $\left(\ref{30}\right)$--$\left(\ref{32}\right)$ and perform the integrations
%
\begin{eqnarray}
\label{64}
\mathcal{O}_i^{(l)}\left( {x,t} \right)=\frac{1}{{{{\left( {2\sqrt {\pi \kappa } } \right)}^n}}}\int\limits_0^t {\int\limits_{{\mathbb{R}^n}} {\frac{{q_i^{(l)}\left( {y,t} \right){e^{ - \sum\limits_{k = 1}^n {\frac{{{{\left( {{x_k} - {y_k}} \right)}^2}}}{{4\kappa \left( {t - \tau } \right)}}} }}}}{{{{\left( {t - \tau } \right)}^{\frac{n}{2}}}}}\prod\limits_{j = 1}^n {d{y_j}} } } d\tau,\quad l=2,3,... 
\end{eqnarray}
$\bm{\left(iii\right)}$ Obtain $v_i^{(l)}$ from the formula
%
\begin{eqnarray}
\label{65}
v_i^{(1)}\left( {x,t} \right) &=& v_i^{(0)}\left( x \right){\mathcal{T} _1}\left( t \right)\nonumber\\
v_i^{(l)}\left( {x,t} \right) &=& v_i^{(l - 1)}\left( {x,t} \right) - \mathcal{O}_i^{(l-1)}\left( {x,t} \right),\quad l=2,3,...
\end{eqnarray}
The generating function $\mathcal{O}_i^{(l - 1)}\left( {x,t} \right)$ takes the following form:
%
\begin{eqnarray}
\label{66}
\mathcal{O}_i^{(l)}\left( {x,t} \right) = \sum\limits_{j = 1}^{m} {\chi  _{ij}^{(l)}\!\left( x \right)\mathcal{T}_{j+m}\!\left( t \right)} 
\end{eqnarray}
Where $m$ is the number of additional terms generated in the $l^{th}$ sequence, $\chi  _{ij}^{(l)}\left( x \right)$ is a spatially periodic function and $\mathcal{T}_{j+m}\!\left( t \right)$ is a product of $\left( {{\pi}\kappa } \right)$ with a negative exponent and the sum of a finite series involving exponentially decaying functions in time of the form ${e^{ - \varepsilon {\pi ^2}\kappa t}}$ and ${\left( {\kappa t} \right)^\sigma }{e^{ - \varepsilon {\pi ^2}\kappa t}}$, $\varepsilon  \gg \sigma $. Where $\varepsilon$ and $\sigma$ are integers. It is apparent that, for a given positive coefficient of kinematical viscosity, $\kappa$, it is the form of the time only function $\mathcal{T}_{j+m}\!\left( t \right)$ that determines the extent of the time in which the solutions accelerate away towards a higher value (not a singularity) and quickly recedes, a phenomenon known as \protect{\itshape{blowup time}} (See example in the Appendices).\\\\
Pressure is given by
%
\begin{eqnarray}
\label{67}
p^{(l)}\left( {x,t} \right) =   \rho \int {\hbar _i^{(l - 1)}d{x_i}}  + C\left( {{x_1},{x_2},...,{x_{i - 1}},{x_{i + 1}},...,{x_n}} \right)
\end{eqnarray}
where ${C}$ is the integration constant.\\\\
We conclude that when $\mathcal{O}_i^{(l)}\left( {x,t} \right)$ becomes vanishingly small the velocity vector field $v_i^{(l)}\left( {x,t} \right)$ at sequence $l$ satisfies the Navier-Stokes system of equations $\left(\ref{1}\right)-\left(\ref{4}\right)$. At very low Reynolds numbers, $\mathcal{R}e \ll 1$, the viscous forces dominate over the inertia forces.  Thus, the latter may be neglected in the Navier-Stokes equations. It is implicit form $\left(\ref{66}\right)$ that at small $\mathcal{R}e$, only a few sequences will be required to make $\mathcal{O}_i^{(l)}\left( {x,t} \right)$ vanish. As the $\mathcal{R}e$ increases, more and more sequences will be required before $\mathcal{O}_i^{(l)}\left( {x,t} \right)$ would become vanishingly small. Nonetheless, the theory holds for all $\kappa  > 0$.\\\\
The sequence by sequence process of deriving analytic expressions of $v_i^{(l)}\left( {x,t} \right)$, though straightforward, are exhaustively lengthy and time consuming. Mathematical tools such as \protect{\itshape{Mathematica}} or  MATLAB may be used to perform symbolic manipulations to express $q_i^{(l)}\left( {x,t} \right)$ in a form integrable by use of the identities $\left(\ref{30}\right)$--$\left(\ref{32}\right)$.\\\\
In the next section, for a given solenoidal initial velocity vector field in $\mathbb{R}^3$, we derive three sequences of velocity vector fields to show that a consistent pattern, subjecting the \protect{\itshape{blowup time}}, develops with increasing $\mathcal{R}e$. The expressions derived for velocity and pressure are smooth and satisfy $\left(\ref{1}\right)-\left(\ref{4}\right)$ in the applicable range of the $\mathcal{R}e$.
\vspace{-0.1in}
\section{An illustrative example in $\mathbb{R}^3$}
We choose a spatially symmetric initial condition and derive three sequences of the velocity field.
%
\begin{eqnarray}
\label{68}
v_1^{(0)}  = \sin \left( {\pi x_1 } \right)\sin \left( {\pi x_3 } \right) + \cos \left( {\pi x_1 } \right)\cos \left( {\pi x_2 } \right)
\end{eqnarray}
%
\begin{eqnarray}
\label{69}
v_2^{(0)}  = \sin \left( {\pi x_2 } \right)\sin \left( {\pi x_1 } \right) + \cos \left( {\pi x_2 } \right)\cos \left( {\pi x_3 } \right)
\end{eqnarray}
%
\begin{eqnarray}
\label{70}
v_3^{(0)}  = \sin \left( {\pi x_3 } \right)\sin \left( {\pi x_2 } \right) + \cos \left( {\pi x_3 } \right)\cos \left( {\pi x_1 } \right)
\end{eqnarray}
$v_i^{(0)}$, $i=1,2,3$ is a smooth spatially periodic vector field satisfying $\left(\ref{2}\right)$. We take $f_i\left( {x,t} \right)$ to be identically zero.\\\\
{\protect{\bfseries{Sequence: 1}}}\\\\
Evaluating the integrals in  $\left(\ref{29}\right)$ by use of the integral identities $\left(\ref{30}\right)$--$\left(\ref{32}\right)$ we find $\xi_1=2$. 
%
\begin{eqnarray}
\label{71}
v_i^{(1)} = v_i^{(0)}\left( x \right){\mathcal{T} _1}\left( t \right)
\end{eqnarray}
where
%
\begin{eqnarray}
\label{72}
{\mathcal{T} _1}\left( t \right) = {e^{ - {2}{\pi ^2}\kappa t}}
\end{eqnarray}
From $\left(\ref{36}\right)$ we get
%
\begin{eqnarray}
\label{73}
\mathcal{U}_1^{(0)}&=&g_1^{(0)}  - \frac{1}
{{4\pi }}\int\limits_{\mathbb{R}^3 } {\frac{{\left( {x_1  - y_1 } \right)\sum\limits_{k = 1}^3 {\frac{{\partial g_k^{(0)} \left( {y ,t} \right)}}
{{\partial y_k }}} }}
{{\left\{ {\mathcal{P}_3\left( {x,y} \right)} \right\}^{\frac{3}
{2}} }}} \prod\limits_{j = 1}^3 {dy_j }\nonumber\\
&=&\frac{{2\pi }}
{3}\cos \left( {2\pi x_1 } \right)\cos \left( {\pi x_2 } \right)\sin \left( {\pi x_3 } \right) - \frac{{2\pi }}
{3}\sin \left( {2\pi x_2 } \right)\cos \left( {\pi x_1 } \right)\cos \left( {\pi x_3 } \right) +\nonumber\\
&+&\frac{{2\pi }}
{3}\sin \left( {2\pi x_3 } \right)\sin \left( {\pi x_1 } \right)\sin \left( {\pi x_2 } \right)
\end{eqnarray}
%
\begin{eqnarray}
\label{74}
\mathcal{U}_2^{(0)} &  = & g_2^{(0)}-\frac{1}
{{4\pi }}\int\limits_{\mathbb{R}^3 } {\frac{{\left( {x_2  - y_2 } \right)\sum\limits_{k = 1}^3 {\frac{{\partial g_k^{(0)} \left( {y ,t} \right)}}
{{\partial y_k }}} }}
{{\left\{ {\mathcal{P}_3\left( {x,y} \right)} \right\}^{\frac{3}
{2}} }}} \prod\limits_{j = 1}^3 {dy_j }\nonumber\\
&=&\frac{{2\pi }}
{3}\cos \left( {2\pi x_2 } \right)\sin \left( {\pi x_1 } \right)\cos \left( {\pi x_3 } \right) - \frac{{2\pi }}
{3}\sin \left( {2\pi x_3 } \right)\cos \left( {\pi x_1 } \right)\cos \left( {\pi x_2 } \right)
 +\nonumber\\
&+&\frac{{2\pi }}
{3}\sin \left( {2\pi x_1 } \right)\sin \left( {\pi x_3 } \right)\sin \left( {\pi x_2 } \right)
\end{eqnarray}
%
\begin{eqnarray}
\label{75}
\mathcal{U}_3^{(0)} &  = &g_3^{(0)}-\frac{1}
{{4\pi }}\int\limits_{\mathbb{R}^3 } {\frac{{\left( {x_3  - y_3 } \right)\sum\limits_{k = 1}^3 {\frac{{\partial g_k^{(0)} \left( {y,t} \right)}}
{{\partial y_k}}} }}
{{\left\{ {\mathcal{P}_3\left( {x,y } \right)} \right\}^{\frac{3}
{2}} }}} \prod\limits_{j = 1}^3 {dy_j }\nonumber\\
&=&
\frac{{2\pi }}
{3}\cos \left( {2\pi x_3 } \right)\cos \left( {\pi x_1 } \right)\sin \left( {\pi x_2 } \right) - \frac{{2\pi }}
{3}\sin \left( {2\pi x_1 } \right)\cos \left( {\pi x_2 } \right)\cos \left( {\pi x_3 } \right)
 +\nonumber\\
&+&\frac{{2\pi }}
{3}\sin \left( {2\pi x_2 } \right)\sin \left( {\pi x_1 } \right)\sin \left( {\pi x_3 } \right)
\end{eqnarray}
Where $g_i^{(0)}$, $i=1,2,3$, is given by $\left(\ref{37}\right)$. In evaluating the integrals in $\left(\ref{73}\right)$, $\left(\ref{74}\right)$ and $\left(\ref{75}\right)$ we have used the following identities: 
\begin{eqnarray}
\label{76}
\begin{gathered}
  \int\limits_0^\infty  {\frac{{u\sin \left( {\alpha u} \right)}}{{\sqrt {{{\left( {{\beta ^2} + {u^2}} \right)}^3}} }}} du = \alpha {K_0}\left( {\alpha \beta } \right) \hfill \\
  \int\limits_0^\infty  {{K_0}\left( {\alpha \sqrt {{u^2} + {\beta ^2}} } \right)\cos \left( {z u} \right)du}  = \frac{{\pi {e^{ - \beta \sqrt {{\alpha ^2} + {z ^2}} }}}}{{2\sqrt {{\alpha ^2} + {z ^2}} }} \hfill \\ 
\end{gathered} 
\end{eqnarray}
Where $K_\nu  \left( {u } \right)$ is the modified Bessel function of the second kind of order $\nu$.\\\\
{\protect{\bfseries{Sequence: 2}}}
%
\begin{eqnarray}
\label{77}
\mathcal{U}_i^{(1)} = g_i^{(1)} - \hbar _i^{(1)} = \mathcal{U}_i^{(0)}{e^{ - 4{\pi ^2}\kappa t}}
\end{eqnarray}
Evaluating the integrals in  $\left(\ref{44}\right)$ by use of the integral identities $\left(\ref{30}\right)$--$\left(\ref{32}\right)$ we find $\xi_2=6$. 
%
\begin{eqnarray}
\label{78}
v_i^{(2)} = v_i^{(1)} - \mathcal{U}_i^{(0)}{\mathcal{T} _2}\left( t \right)
\end{eqnarray}
%
Where
%
\begin{eqnarray}
\label{79}
{\mathcal{T} _2}\left( t \right) = \frac{{\left( {{e^{ - 4{\pi ^2}\kappa t}} - {e^{ - 6{\pi ^2}\kappa t}}} \right)}}{{2{\pi ^2}\kappa }}
\end{eqnarray}
{\protect{\bfseries{Sequence: 3}}}
%
\begin{eqnarray}
\label{80}
\mathcal{U}_i^{(2)} = g_i^{(2)} - \hbar _i^{(2)} 
\end{eqnarray}
where $g_i^{(2)}$ is given by $\left(\ref{48}\right)$, which is
%
\begin{eqnarray}
\label{81}
g_i^{(2)} = \sum\limits_{j = 1}^n {v_j^2} \frac{{\partial v_i^{(2)}}}{{\partial {x_j}}} = g_i^{(1)} + q_i^{(2)}
\end{eqnarray}
We note that, in this particular example,  ${\mathbf{div}}{\kern 1pt} {\kern 1pt} \bm{q}^{(2)} = \sum\limits_{i = 1}^n {\frac{{\partial q_i^{(2)}}}{{\partial {x_i}}}}  = 0$. Therefore, 
\begin{center}
$\begin{gathered}
  \beta _{ik}^{(2)}\left( x \right) = 0,\,\,k = 1,2 \hfill \\
  \mathcal{G}_{ik}^{(2)}\left( x \right) = \alpha _{ik}^{(2)}\left( x \right),\,\,k = 1,2 \hfill \\
  \hbar _i^{(2)} = \hbar _i^{(1)} \hfill \\ 
\end{gathered}$
\end{center}
Equation $\left(\ref{55}\right)$ simplifies to
%
\begin{eqnarray}
\label{82}
q_i^{(2)} = \alpha _{i1}^{(2)}\left(x\right){\mathcal{T} _1}\left( t \right){\mathcal{T} _2}\left( t \right) + \alpha _{i2}^{(2)}\left(x\right){\mathcal{T} _2^2\left( t \right)}
\end{eqnarray}
Where
\begin{eqnarray}
\label{83}
{\mathcal{T} _1}\left( t \right){\mathcal{T} _2}\left( t \right) = \frac{{\left( {{e^{ - 6{\pi ^2}\kappa t}} - {e^{ - 8{\pi ^2}\kappa t}}} \right)}}{{2{\pi ^2}\kappa }}
\end{eqnarray}
and
%
\begin{eqnarray}
\label{84}
{\mathcal{T} _2^2\left( t \right)} = \frac{{\left( {{e^{ - 8{\pi ^2}\kappa t}} - 2{e^{ - 10{\pi ^2}\kappa t}} + {e^{ - 12{\pi ^2}\kappa t}}} \right)}}{{4{\pi ^4}{\kappa ^2}}}
\end{eqnarray}
The coefficients $\alpha _{ik}^{(2)}\left(x\right)$, $k=1,2$ are obtained from $\left(\ref{50}\right)$ and $\left(\ref{51}\right)$ and are given, in a form integrable by use of the identities $\left(\ref{30}\right)$--$\left(\ref{32}\right)$, in the \protect{\itshape{Appendices}}. We have\\\\ 
%
\begin{eqnarray}
\label{85}
\mathcal{U}_i^{(2)}= g_i^{(1)} + q_i^{(2)} - \hbar _i^{(1)} = \mathcal{U}_i^{(1)}+ q_i^{(2)}
\end{eqnarray}
and 
%
\begin{eqnarray}
\label{86}
v_i^{(3)} = v_i^{(2)} - \mathcal{O}_i^{(2)}
\end{eqnarray}
%
Where
%
\begin{eqnarray}
\label{87}
\mathcal{O}_i^{(2)}\left( {x,t} \right) &=& \frac{1}{{{{\left( {2\sqrt {\pi \kappa } } \right)}^n}}}\int\limits_0^t {\int\limits_{{\mathbb{R}^n}} {\frac{{q_i^{(2)}\left( {y,t} \right){e^{ - \sum\limits_{k = 1}^n {\frac{{{{\left( {{x_k} - {y_k}} \right)}^2}}}{{4\kappa \left( {t - \tau } \right)}}} }}}}{{{{\left( {t - \tau } \right)}^{\frac{n}{2}}}}}\prod\limits_{j = 1}^n {d{y_j}} } } d\tau\nonumber\\
&=&\sum\limits_{j = 1}^{11} {\chi  _{ij}^{(2)}\left( x \right){\mathcal{T} _{j + 2}}\left( t \right)} 
\end{eqnarray}
The coefficients $\chi  _{ij}^{(2)}\left( x \right)$ and $\mathcal{T} _j\left( {t} \right)$,  $i=1,2,3$, $j=1,11$ are given, in a form integrable by use of the identities $\left(\ref{30}\right)$--$\left(\ref{32}\right)$, in the \protect{\itshape{Appendices}}.\\\\ 

Pressure at the end of the third sequence is obtained from
%
\begin{eqnarray}
\label{89}
p^{(3)} =   \rho \int {\hbar _1^{(2)}d{x_1}}  + C\left( {{x_2},{x_3}} \right)
\end{eqnarray}
Figures 1-12 shows plots of $g_1^{(k)}$ and $v_1^{(k)}$, $k=0,1,2,3$ versus log of time for progressively increasing $\mathcal{R}e$ at an arbitrarily chosen point $\left( {{x_1},{x_2},{x_3}} \right) \equiv \left( {\frac{{2\pi }}{3},\frac{{3\pi }}{5},\frac{{4\pi }}{9}} \right)$\footnote{Computations have been performed with dimensionless variables: $\kappa$ is replaced by $\frac{1}{{\mathcal{R}e}}$}. The results, for this particular example, is self-explanatory. $\mathcal{R}e < 1$ requires only one sequence and, as expected, for ``Creeping flow" where the nonlinear inertial convective term is negligible in comparison to the rest of the terms, the solution, $v_i^{(1)}\left( {x,t} \right)$, is given by the linear homogeneous diffusion equation.  $\mathcal{R}e < 15$ and $\mathcal{R}e  < 20$ may require two and three sequences respectively. $\mathcal{R}e >20$ will require more sequences to be computed.\\\\
The purpose of these simple illustrations is not to define ranges of $\mathcal{R}e$ and determine the corresponding number of sequences required to arrive at a solution, but rather to simply demonstrate that a pattern of behavior, including that of the \protect{\itshape{blowup time}}, that develops with increasing $\mathcal{R}e$.\\\\\\
%
%
\begin{figure}
\begin{center}
\begin{subfigure}[]{0.48\textwidth}
\centering
\resizebox{2.50in}{!}{
 \rotatebox{0}{
\includegraphics{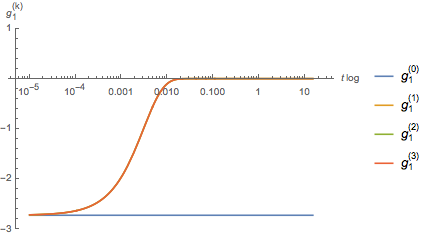}
}}
\caption{$g_i^{(k)}$, $ k=0,1,2,3$ vs $\log \left(t\right)$}
\end{subfigure}
\begin{subfigure}[]{0.48\textwidth}
\centering
\resizebox{2.40in}{!}{
 \rotatebox{0}{
\includegraphics{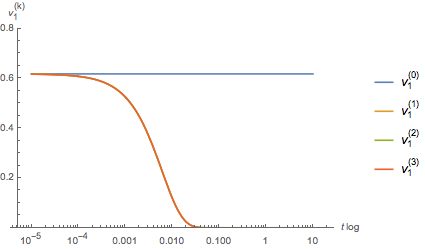}
}}
\caption{$v_i^{(k)}$, $ k=0,1,2,3$ vs $\log \left(t\right)$}
\end{subfigure}
\end{center}
\caption{Creeping flow --- $\mathcal{R}e \leqslant 1$, ($\mathcal{R}e=0.06$)}
\end{figure}
%
%
\begin{figure}
\begin{center}
\begin{subfigure}[]{0.48\textwidth}
\centering
\resizebox{2.50in}{!}{
 \rotatebox{0}{
\includegraphics{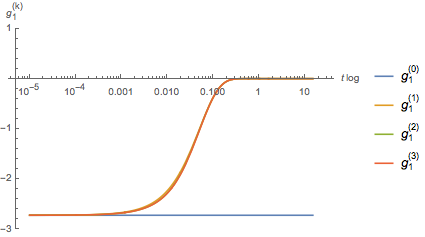}
}}
\caption{$g_i^{(k)}$, $ k=0,1,2,3$ vs $\log \left(t\right)$}
\end{subfigure}
\begin{subfigure}[]{0.48\textwidth}
\centering
\resizebox{2.40in}{!}{
 \rotatebox{0}{
\includegraphics{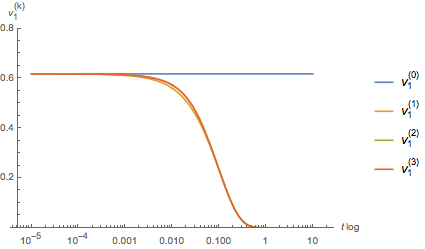}
}}
\caption{$v_i^{(k)}$, $ k=0,1,2,3$ vs $\log \left(t\right)$}
\end{subfigure}
\end{center}
\caption{$\mathcal{R}e =1$}
\end{figure}
%
%
\begin{figure}
\begin{center}
\begin{subfigure}[]{0.48\textwidth}
\centering
\resizebox{2.50in}{!}{
 \rotatebox{0}{
\includegraphics{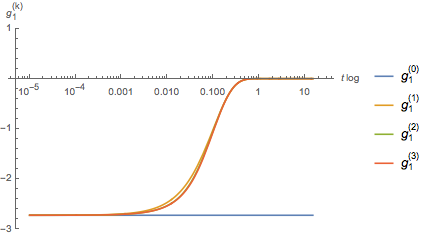}
}}
\caption{$g_i^{(k)}$, $ k=0,1,2,3$ vs $\log \left(t\right)$}
\end{subfigure}
\begin{subfigure}[]{0.48\textwidth}
\centering
\resizebox{2.40in}{!}{
 \rotatebox{0}{
\includegraphics{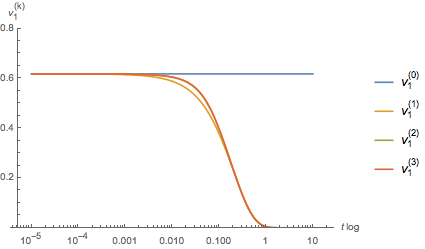}
}}
\caption{$v_i^{(k)}$, $ k=0,1,2,3$ vs $\log \left(t\right)$}
\end{subfigure}
\end{center}
\caption{$\mathcal{R}e =2$}
\end{figure}
%
%
\begin{figure}
\begin{center}
\begin{subfigure}[]{0.48\textwidth}
\centering
\resizebox{2.50in}{!}{
 \rotatebox{0}{
\includegraphics{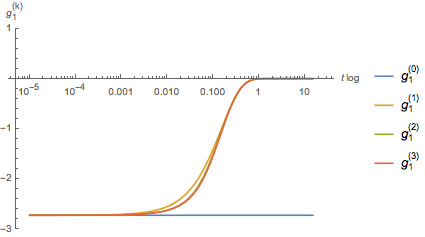}
}}
\caption{$g_i^{(k)}$, $ k=0,1,2,3$ vs $\log \left(t\right)$}
\end{subfigure}
\begin{subfigure}[]{0.48\textwidth}
\centering
\resizebox{2.40in}{!}{
 \rotatebox{0}{
\includegraphics{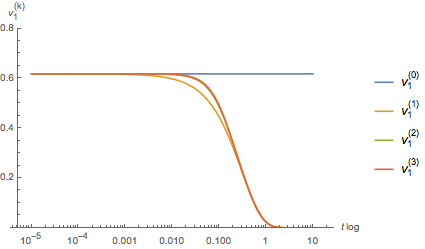}
}}
\caption{$v_i^{(k)}$, $ k=0,1,2,3$ vs $\log \left(t\right)$}
\end{subfigure}
\end{center}
\caption{$\mathcal{R}e =3$}
\end{figure}
%
%
\begin{figure}
\begin{center}
\begin{subfigure}[]{0.48\textwidth}
\centering
\resizebox{2.50in}{!}{
 \rotatebox{0}{
\includegraphics{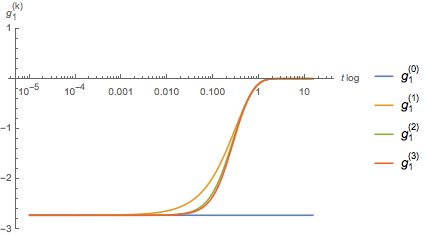}
}}
\caption{$g_i^{(k)}$, $ k=0,1,2,3$ vs $\log \left(t\right)$}
\end{subfigure}
\begin{subfigure}[]{0.48\textwidth}
\centering
\resizebox{2.40in}{!}{
 \rotatebox{0}{
\includegraphics{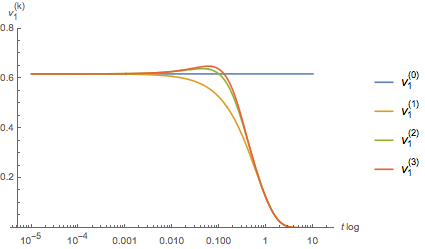}
}}
\caption{$v_i^{(k)}$, $ k=0,1,2,3$ vs $\log \left(t\right)$}
\end{subfigure}
\end{center}
\caption{$\mathcal{R}e =6$}
\end{figure}
%
%
\begin{figure}
\begin{center}
\begin{subfigure}[]{0.48\textwidth}
\centering
\resizebox{2.50in}{!}{
 \rotatebox{0}{
\includegraphics{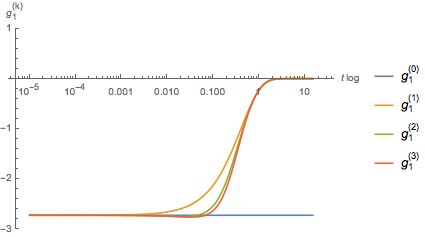}
}}
\caption{$g_i^{(k)}$, $ k=0,1,2,3$ vs $\log \left(t\right)$}
\end{subfigure}
\begin{subfigure}[]{0.48\textwidth}
\centering
\resizebox{2.40in}{!}{
 \rotatebox{0}{
\includegraphics{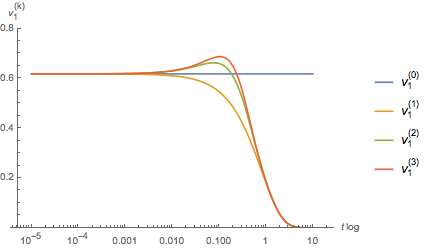}
}}
\caption{$v_i^{(k)}$, $ k=0,1,2,3$ vs $\log \left(t\right)$}
\end{subfigure}
\end{center}
\caption{$\mathcal{R}e =8$}
\end{figure}
%
%
\begin{figure}
\begin{center}
\begin{subfigure}[]{0.48\textwidth}
\centering
\resizebox{2.50in}{!}{
 \rotatebox{0}{
\includegraphics{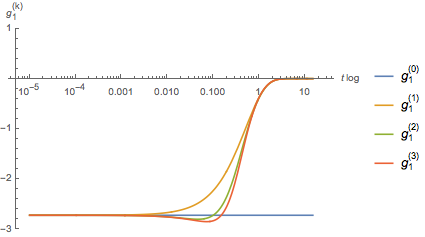}
}}
\caption{$g_i^{(k)}$, $ k=0,1,2,3$ vs $\log \left(t\right)$}
\end{subfigure}
\begin{subfigure}[]{0.48\textwidth}
\centering
\resizebox{2.40in}{!}{
 \rotatebox{0}{
\includegraphics{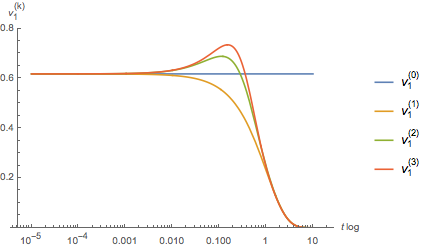}
}}
\caption{$v_i^{(k)}$, $ k=0,1,2,3$ vs $\log \left(t\right)$}
\end{subfigure}
\end{center}
\caption{$\mathcal{R}e =10$}
\end{figure}
%
%
\begin{figure}
\begin{center}
\begin{subfigure}[]{0.48\textwidth}
\centering
\resizebox{2.50in}{!}{
 \rotatebox{0}{
\includegraphics{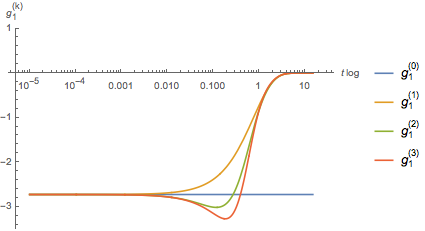}
}}
\caption{$g_i^{(k)}$, $ k=0,1,2,3$ vs $\log \left(t\right)$}
\end{subfigure}
\begin{subfigure}[]{0.48\textwidth}
\centering
\resizebox{2.40in}{!}{
 \rotatebox{0}{
\includegraphics{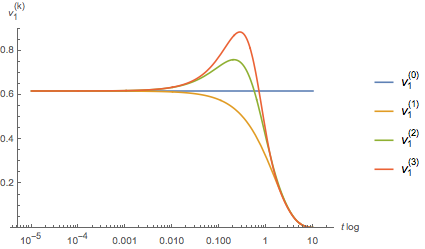}
}}
\caption{$v_i^{(k)}$, $ k=0,1,2,3$ vs $\log \left(t\right)$}
\end{subfigure}
\end{center}
\caption{$\mathcal{R}e =15$}
\end{figure}
%
%
\begin{figure}
\begin{center}
\begin{subfigure}[]{0.48\textwidth}
\centering
\resizebox{2.50in}{!}{
 \rotatebox{0}{
\includegraphics{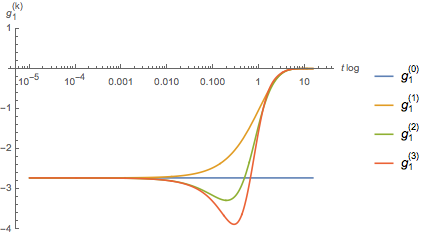}
}}
\caption{$g_i^{(k)}$, $ k=0,1,2,3$ vs $\log \left(t\right)$}
\end{subfigure}
\begin{subfigure}[]{0.48\textwidth}
\centering
\resizebox{2.40in}{!}{
 \rotatebox{0}{
\includegraphics{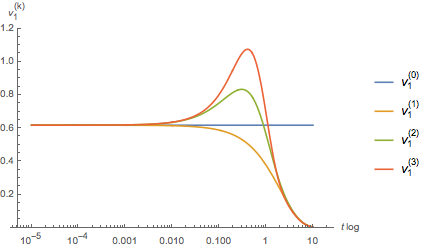}
}}
\caption{$v_i^{(k)}$, $ k=0,1,2,3$ vs $\log \left(t\right)$}
\end{subfigure}
\end{center}
\caption{$\mathcal{R}e =20$}
\end{figure}
%
%
\begin{figure}
\begin{center}
\begin{subfigure}[]{0.48\textwidth}
\centering
\resizebox{2.50in}{!}{
 \rotatebox{0}{
\includegraphics{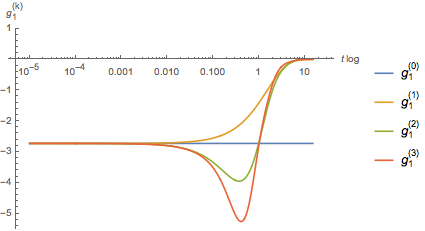}
}}
\caption{$g_i^{(k)}$, $ k=0,1,2,3$ vs $\log \left(t\right)$}
\end{subfigure}
\begin{subfigure}[]{0.48\textwidth}
\centering
\resizebox{2.40in}{!}{
 \rotatebox{0}{
\includegraphics{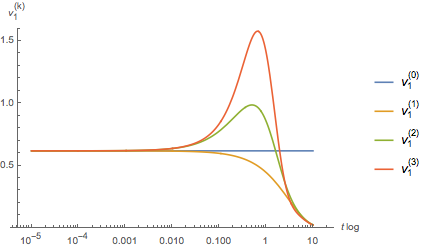}
}}
\caption{$v_i^{(k)}$, $ k=0,1,2,3$ vs $\log \left(t\right)$}
\end{subfigure}
\end{center}
\caption{$\mathcal{R}e =30$}
\end{figure}
%
%
\begin{figure}
\begin{center}
\begin{subfigure}[]{0.48\textwidth}
\centering
\resizebox{2.50in}{!}{
 \rotatebox{0}{
\includegraphics{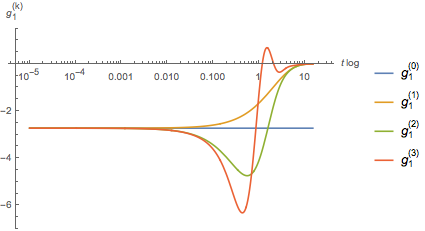}
}}
\caption{$g_i^{(k)}$, $ k=0,1,2,3$ vs $\log \left(t\right)$}
\end{subfigure}
\begin{subfigure}[]{0.48\textwidth}
\centering
\resizebox{2.40in}{!}{
 \rotatebox{0}{
\includegraphics{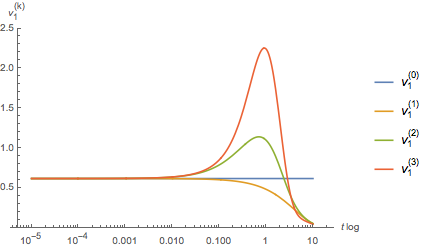}
}}
\caption{$v_i^{(k)}$, $ k=0,1,2,3$ vs $\log \left(t\right)$}
\end{subfigure}
\end{center}
\caption{$\mathcal{R}e =40$}
\end{figure}
%
%
\begin{figure}
\begin{center}
\begin{subfigure}[]{0.48\textwidth}
\centering
\resizebox{2.50in}{!}{
 \rotatebox{0}{
\includegraphics{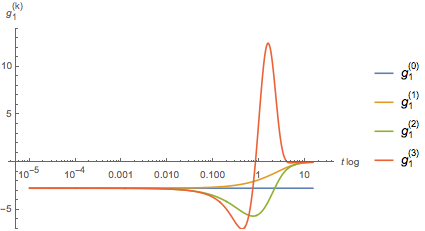}
}}
\caption{$g_i^{(k)}$, $ k=0,1,2,3$ vs $\log \left(t\right)$}
\end{subfigure}
\begin{subfigure}[]{0.48\textwidth}
\centering
\resizebox{2.40in}{!}{
 \rotatebox{0}{
\includegraphics{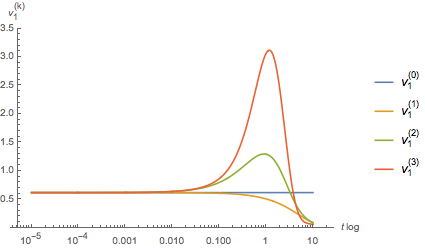}
}}
\caption{$v_i^{(k)}$, $ k=0,1,2,3$ vs $\log \left(t\right)$}
\end{subfigure}
\end{center}
\caption{$\mathcal{R}e =50$}
\end{figure}
\newpage
\section{Concluding Remarks} 
We have presented a solution method that can be used to derive solutions of the non-stationary Navier-Stokes equations $\left(\ref{1}\right)-\left(\ref{4}\right)$ for incompressible viscous fluids in $\mathbb{R}^n$. The essence of the solution method can be summarized as follows:\\\\
$\bm{\left(i\right)}$ We recast the Navier-Stokes equation for velocity in terms of three distinct terms associated, respectively, with the linear viscous force, the externally applied force and  the inertial force, given by $\left(\ref{13}\right)$.\\\\
$\bm{\left(ii\right)}$ We observe that, for a given solenoidal initial velocity vector field, the nonlinearity, expressed in terms of the components of ${\mathcal{U}_i}\left( x,t \right)$, \protect{\itshape{instantaneously}} spirals from zero to a plateau through a sequence of linear diffusion processes in accordance with $\left(\ref{13}\right)$. An analytical expression for the plateaued ${\mathcal{U}_i}\left( x,t \right)$ is given by $\left(\ref{63}\right)$. The plateaued ${\mathcal{U}_i}\left( x,t \right)$, when substituted into $\left(\ref{13}\right)$ results in a linear inhomogeneous partial differential equation, which is solved analytically.\\\\
$\bm{\left(iii\right)}$ As the $\mathcal{R}e$ increases, more and more sequences will be required to arrive at a complete solution. The closed-form analytic solution is composed of a  finite series. We show that the number of terms required in the finite series is dependent on the $\mathcal{R}e$. If the $\mathcal{R}e$ is small, the solution will contain fewer terms. As the $\mathcal{R}e$ increases, the number of terms required to complete the closed form solution will also increase correspondingly.\\\\
$\bm{\left(iv\right)}$ We show that for a given positive coefficient of kinematical viscosity, $\kappa$, it is the form of the time only function $T_{j+m}\!\left( t \right)$ that determines the extent of the time in which the solution accelerate towards a higher value (not a singularity), a phenomenon known as \protect{\itshape{blowup time}}. The solution, once past this higher value, quickly recedes.\\\\
$\bm{\left(v\right)}$ The pressure field is given by the solution of the Poisson equation, $\left(\ref{8}\right)$.\\\\
$\bm{\left(vi\right)}$ The solution presented in $\mathbb{R}^3$ for velocity and pressure are smooth and satisfy $\left(\ref{1}\right)-\left(\ref{4}\right)$.\\
\bibliographystyle{chicago}
\bibliography{NS}
\section{Appendices}
\subsection{\protect{\bfseries{The coefficients $\alpha _{ik}^{(2)}\left(x\right)$, $i=1,2,3,\,\,k=1,2$ in $\left(\ref{82}\right)$}}}
%
\begin{eqnarray}
\label{90}
\alpha _{11}^{(2)}\left( x \right)& =& \frac{{{\pi ^2}}}{3}\sin \left( {\pi {x_1}} \right)\sin \left( {\pi {x_3}} \right) + \frac{{{\pi ^2}}}{3}\cos \left( {\pi {x_1}} \right)\cos \left( {\pi {x_2}} \right) + \nonumber\\
&+& \frac{{{\pi ^2}}}{3}\cos \left( {\pi {x_1}} \right)\cos \left( {\pi {x_2}} \right)\cos \left( {2\pi {x_3}} \right) - \frac{{{\pi ^2}}}{3}\sin \left( {\pi {x_1}} \right)\cos \left( {2\pi {x_2}} \right)\sin \left( {\pi {x_3}} \right) - \nonumber\\
& -& \frac{{5{\pi ^2}}}{3}\sin \left( {2\pi {x_1}} \right)\sin \left( {\pi {x_2}} \right)\cos \left( {\pi {x_3}} \right) - \frac{{2{\pi ^2}}}{3}\sin \left( {2\pi {x_2}} \right)\sin \left( {2\pi {x_3}} \right) +\nonumber\\
&+&\frac{{5{\pi ^2}}}{6}\sin \left( {3\pi {x_1}} \right)\sin \left( {\pi {x_3}} \right) - \frac{{5{\pi ^2}}}{6}\cos \left( {3\pi {x_1}} \right)\cos \left( {\pi {x_2}} \right) + \nonumber\\
&+& \frac{{{\pi ^2}}}{2}\cos \left( {\pi {x_1}} \right)\cos \left( {3\pi {x_2}} \right) - \frac{{{\pi ^2}}}{2}\sin \left( {\pi {x_1}} \right)\sin \left( {3\pi {x_3}} \right) + \nonumber\\
&+& \frac{{{\pi ^2}}}{6}\cos \left( {3\pi {x_1}} \right)\cos \left( {\pi {x_2}} \right)\cos \left( {2\pi {x_3}} \right) + \frac{{{\pi ^2}}}{6}\cos \left( {\pi {x_1}} \right)\cos \left( {3\pi {x_2}} \right)\cos \left( {2\pi {x_3}} \right) + \nonumber\\
&+& \frac{{{\pi ^2}}}{6}\sin \left( {\pi {x_1}} \right)\cos \left( {2\pi {x_2}} \right)\sin \left( {3\pi {x_3}} \right) - \frac{{{\pi ^2}}}{6}\sin \left( {2\pi {x_1}} \right)\sin \left( {\pi {x_2}} \right)\cos \left( {3\pi {x_3}} \right) + \nonumber\\
&+& \frac{{{\pi ^2}}}{6}\sin \left( {2\pi {x_1}} \right)\sin \left( {3\pi {x_2}} \right)\cos \left( {\pi {x_3}} \right) + \frac{{{\pi ^2}}}{6}\sin \left( {3\pi {x_1}} \right)\cos \left( {2\pi {x_2}} \right)\sin \left( {\pi {x_3}} \right)\qquad
\end{eqnarray}
%
\begin{eqnarray}
\label{91}
\alpha _{21}^{(2)}\left( x \right)& =& \frac{{{\pi ^2}}}{3}\sin \left( {\pi {x_2}} \right)\sin \left( {\pi {x_1}} \right) + \frac{{{\pi ^2}}}{3}\cos \left( {\pi {x_2}} \right)\cos \left( {\pi {x_3}} \right) + \nonumber\\
&+& \frac{{{\pi ^2}}}{3}\cos \left( {\pi {x_2}} \right)\cos \left( {\pi {x_3}} \right)\cos \left( {2\pi {x_1}} \right) - \frac{{{\pi ^2}}}{3}\sin \left( {\pi {x_2}} \right)\cos \left( {2\pi {x_3}} \right)\sin \left( {\pi {x_1}} \right) - \nonumber\\
& -& \frac{{5{\pi ^2}}}{3}\sin \left( {2\pi {x_2}} \right)\sin \left( {\pi {x_3}} \right)\cos \left( {\pi {x_1}} \right) - \frac{{2{\pi ^2}}}{3}\sin \left( {2\pi {x_3}} \right)\sin \left( {2\pi {x_1}} \right) +\nonumber\\
&+&\frac{{5{\pi ^2}}}{6}\sin \left( {3\pi {x_2}} \right)\sin \left( {\pi {x_1}} \right) - \frac{{5{\pi ^2}}}{6}\cos \left( {3\pi {x_2}} \right)\cos \left( {\pi {x_3}} \right) + \nonumber\\
&+& \frac{{{\pi ^2}}}{2}\cos \left( {\pi {x_2}} \right)\cos \left( {3\pi {x_3}} \right) - \frac{{{\pi ^2}}}{2}\sin \left( {\pi {x_2}} \right)\sin \left( {3\pi {x_1}} \right) + \nonumber\\
&+& \frac{{{\pi ^2}}}{6}\cos \left( {3\pi {x_2}} \right)\cos \left( {\pi {x_3}} \right)\cos \left( {2\pi {x_1}} \right) + \frac{{{\pi ^2}}}{6}\cos \left( {\pi {x_2}} \right)\cos \left( {3\pi {x_3}} \right)\cos \left( {2\pi {x_1}} \right) + \nonumber\\
&+& \frac{{{\pi ^2}}}{6}\sin \left( {\pi {x_2}} \right)\cos \left( {2\pi {x_3}} \right)\sin \left( {3\pi {x_1}} \right) - \frac{{{\pi ^2}}}{6}\sin \left( {2\pi {x_2}} \right)\sin \left( {\pi {x_3}} \right)\cos \left( {3\pi {x_1}} \right) + \nonumber\\
&+& \frac{{{\pi ^2}}}{6}\sin \left( {2\pi {x_2}} \right)\sin \left( {3\pi {x_3}} \right)\cos \left( {\pi {x_1}} \right) + \frac{{{\pi ^2}}}{6}\sin \left( {3\pi {x_2}} \right)\cos \left( {2\pi {x_3}} \right)\sin \left( {\pi {x_1}} \right)\qquad
\end{eqnarray}
%
\begin{eqnarray}
\label{92}
\alpha _{31}^{(2)}\left( x \right)& =& \frac{{{\pi ^2}}}{3}\sin \left( {\pi {x_3}} \right)\sin \left( {\pi {x_2}} \right) + \frac{{{\pi ^2}}}{3}\cos \left( {\pi {x_3}} \right)\cos \left( {\pi {x_1}} \right) + \nonumber\\
&+& \frac{{{\pi ^2}}}{3}\cos \left( {\pi {x_3}} \right)\cos \left( {\pi {x_1}} \right)\cos \left( {2\pi {x_2}} \right) - \frac{{{\pi ^2}}}{3}\sin \left( {\pi {x_3}} \right)\cos \left( {2\pi {x_1}} \right)\sin \left( {\pi {x_2}} \right) - \nonumber\\
& -& \frac{{5{\pi ^2}}}{3}\sin \left( {2\pi {x_3}} \right)\sin \left( {\pi {x_1}} \right)\cos \left( {\pi {x_2}} \right) - \frac{{2{\pi ^2}}}{3}\sin \left( {2\pi {x_1}} \right)\sin \left( {2\pi {x_2}} \right) +\nonumber\\
&+&\frac{{5{\pi ^2}}}{6}\sin \left( {3\pi {x_3}} \right)\sin \left( {\pi {x_2}} \right) - \frac{{5{\pi ^2}}}{6}\cos \left( {3\pi {x_3}} \right)\cos \left( {\pi {x_1}} \right) + \nonumber\\
&+& \frac{{{\pi ^2}}}{2}\cos \left( {\pi {x_3}} \right)\cos \left( {3\pi {x_1}} \right) - \frac{{{\pi ^2}}}{2}\sin \left( {\pi {x_3}} \right)\sin \left( {3\pi {x_2}} \right) + \nonumber\\
&+& \frac{{{\pi ^2}}}{6}\cos \left( {3\pi {x_3}} \right)\cos \left( {\pi {x_1}} \right)\cos \left( {2\pi {x_2}} \right) + \frac{{{\pi ^2}}}{6}\cos \left( {\pi {x_3}} \right)\cos \left( {3\pi {x_1}} \right)\cos \left( {2\pi {x_2}} \right) + \nonumber\\
&+& \frac{{{\pi ^2}}}{6}\sin \left( {\pi {x_3}} \right)\cos \left( {2\pi {x_1}} \right)\sin \left( {3\pi {x_2}} \right) - \frac{{{\pi ^2}}}{6}\sin \left( {2\pi {x_3}} \right)\sin \left( {\pi {x_1}} \right)\cos \left( {3\pi {x_2}} \right) + \nonumber\\
&+& \frac{{{\pi ^2}}}{6}\sin \left( {2\pi {x_3}} \right)\sin \left( {3\pi {x_1}} \right)\cos \left( {\pi {x_2}} \right) + \frac{{{\pi ^2}}}{6}\sin \left( {3\pi {x_3}} \right)\cos \left( {2\pi {x_1}} \right)\sin \left( {\pi {x_2}} \right)\qquad
\end{eqnarray}
\newpage
%
\begin{eqnarray}
\label{90}
\alpha _{12}^{(2)}\left( x \right)& =&- \frac{{2{\pi ^3}}}{9}\sin \left( {2\pi {x_1}} \right)\cos \left( {2\pi {x_2}} \right) - \frac{{2{\pi ^3}}}{9}\sin \left( {2\pi {x_1}} \right)\cos \left( {2\pi {x_3}} \right)+ \nonumber\\
&+&\frac{{{\pi ^3}}}{3}\cos \left( {3\pi {x_2}} \right)\sin \left( {\pi {x_3}} \right) + \frac{{{\pi ^3}}}{3}\cos \left( {\pi {x_2}} \right)\sin \left( {3\pi {x_3}} \right)+\nonumber\\
&+&\frac{{{\pi ^3}}}{3}\cos \left( {3\pi {x_1}} \right)\sin \left( {2\pi {x_2}} \right)\cos \left( {\pi {x_3}} \right) + \frac{{{\pi ^3}}}{3}\sin \left( {3\pi {x_1}} \right)\sin \left( {\pi {x_2}} \right)\sin \left( {2\pi {x_3}} \right) + \nonumber\\
&+&\frac{{{\pi ^3}}}{9}\cos \left( {2\pi {x_1}} \right)\cos \left( {3\pi {x_2}} \right)\sin \left( {\pi {x_3}} \right) - \frac{{{\pi ^3}}}{9}\cos \left( {\pi {x_1}} \right)\sin \left( {2\pi {x_2}} \right)\cos \left( {3\pi {x_3}} \right) - \nonumber\\
&-&\frac{{{\pi ^3}}}{9}\cos \left( {2\pi {x_1}} \right)\cos \left( {\pi {x_2}} \right)\sin \left( {3\pi {x_3}} \right) - \frac{{{\pi ^3}}}{9}\sin \left( {\pi {x_1}} \right)\sin \left( {3\pi {x_2}} \right)\sin \left( {2\pi {x_3}} \right) - \nonumber\\
&-&\frac{{2{\pi ^3}}}{9}\sin \left( {4\pi {x_1}} \right) + \frac{{{\pi ^3}}}{9}\sin \left( {2\pi {x_1}} \right)\cos \left( {4\pi {x_3}} \right) + \frac{{{\pi ^3}}}{9}\sin \left( {4\pi {x_1}} \right)\cos \left( {2\pi {x_3}} \right) - \nonumber\\
&-&\frac{{{\pi ^3}}}{9}\sin \left( {4\pi {x_1}} \right)\cos \left( {2\pi {x_2}} \right) - \frac{{{\pi ^3}}}{9}\sin \left( {2\pi {x_1}} \right)\cos \left( {4\pi {x_2}} \right) +\nonumber\\
&+&\frac{{2{\pi ^3}}}{9}\cos \left( {2\pi {x_1}} \right)\cos \left( {3\pi {x_2}} \right)\sin \left( {3\pi {x_3}} \right) \qquad
\end{eqnarray}
%
\begin{eqnarray}
\label{90}
\alpha _{22}^{(2)}\left( x \right)& =&- \frac{{2{\pi ^3}}}{9}\sin \left( {2\pi {x_2}} \right)\cos \left( {2\pi {x_3}} \right) - \frac{{2{\pi ^3}}}{9}\sin \left( {2\pi {x_2}} \right)\cos \left( {2\pi {x_1}} \right)+ \nonumber\\
&+&\frac{{{\pi ^3}}}{3}\cos \left( {3\pi {x_3}} \right)\sin \left( {\pi {x_1}} \right) + \frac{{{\pi ^3}}}{3}\cos \left( {\pi {x_3}} \right)\sin \left( {3\pi {x_1}} \right)+\nonumber\\
&+&\frac{{{\pi ^3}}}{3}\cos \left( {3\pi {x_2}} \right)\sin \left( {2\pi {x_3}} \right)\cos \left( {\pi {x_1}} \right) + \frac{{{\pi ^3}}}{3}\sin \left( {3\pi {x_2}} \right)\sin \left( {\pi {x_3}} \right)\sin \left( {2\pi {x_1}} \right) + \nonumber\\
&+&\frac{{{\pi ^3}}}{9}\cos \left( {2\pi {x_2}} \right)\cos \left( {3\pi {x_3}} \right)\sin \left( {\pi {x_1}} \right) - \frac{{{\pi ^3}}}{9}\cos \left( {\pi {x_2}} \right)\sin \left( {2\pi {x_3}} \right)\cos \left( {3\pi {x_1}} \right) - \nonumber\\
&-&\frac{{{\pi ^3}}}{9}\cos \left( {2\pi {x_2}} \right)\cos \left( {\pi {x_3}} \right)\sin \left( {3\pi {x_1}} \right) - \frac{{{\pi ^3}}}{9}\sin \left( {\pi {x_2}} \right)\sin \left( {3\pi {x_3}} \right)\sin \left( {2\pi {x_1}} \right) - \nonumber\\
&-&\frac{{2{\pi ^3}}}{9}\sin \left( {4\pi {x_2}} \right) + \frac{{{\pi ^3}}}{9}\sin \left( {2\pi {x_2}} \right)\cos \left( {4\pi {x_1}} \right) + \frac{{{\pi ^3}}}{9}\sin \left( {4\pi {x_2}} \right)\cos \left( {2\pi {x_1}} \right) - \nonumber\\
&-&\frac{{{\pi ^3}}}{9}\sin \left( {4\pi {x_2}} \right)\cos \left( {2\pi {x_3}} \right) - \frac{{{\pi ^3}}}{9}\sin \left( {2\pi {x_2}} \right)\cos \left( {4\pi {x_3}} \right) +\nonumber\\
&+&\frac{{2{\pi ^3}}}{9}\cos \left( {2\pi {x_2}} \right)\cos \left( {3\pi {x_3}} \right)\sin \left( {3\pi {x_1}} \right) \qquad
\end{eqnarray}
\newpage
\begin{eqnarray}
\label{90}
\alpha _{32}^{(2)}\left( x \right)& =&- \frac{{2{\pi ^3}}}{9}\sin \left( {2\pi {x_3}} \right)\cos \left( {2\pi {x_1}} \right) - \frac{{2{\pi ^3}}}{9}\sin \left( {2\pi {x_3}} \right)\cos \left( {2\pi {x_2}} \right)+ \nonumber\\
&+&\frac{{{\pi ^3}}}{3}\cos \left( {3\pi {x_1}} \right)\sin \left( {\pi {x_2}} \right) + \frac{{{\pi ^3}}}{3}\cos \left( {\pi {x_1}} \right)\sin \left( {3\pi {x_2}} \right)+\nonumber\\
&+&\frac{{{\pi ^3}}}{3}\cos \left( {3\pi {x_3}} \right)\sin \left( {2\pi {x_1}} \right)\cos \left( {\pi {x_2}} \right) + \frac{{{\pi ^3}}}{3}\sin \left( {3\pi {x_3}} \right)\sin \left( {\pi {x_1}} \right)\sin \left( {2\pi {x_2}} \right) + \nonumber\\
&+&\frac{{{\pi ^3}}}{9}\cos \left( {2\pi {x_3}} \right)\cos \left( {3\pi {x_1}} \right)\sin \left( {\pi {x_2}} \right) - \frac{{{\pi ^3}}}{9}\cos \left( {\pi {x_3}} \right)\sin \left( {2\pi {x_1}} \right)\cos \left( {3\pi {x_2}} \right) - \nonumber\\
&-&\frac{{{\pi ^3}}}{9}\cos \left( {2\pi {x_3}} \right)\cos \left( {\pi {x_1}} \right)\sin \left( {3\pi {x_2}} \right) - \frac{{{\pi ^3}}}{9}\sin \left( {\pi {x_3}} \right)\sin \left( {3\pi {x_1}} \right)\sin \left( {2\pi {x_2}} \right) - \nonumber\\
&-&\frac{{2{\pi ^3}}}{9}\sin \left( {4\pi {x_3}} \right) + \frac{{{\pi ^3}}}{9}\sin \left( {2\pi {x_3}} \right)\cos \left( {4\pi {x_2}} \right) + \frac{{{\pi ^3}}}{9}\sin \left( {4\pi {x_3}} \right)\cos \left( {2\pi {x_2}} \right) - \nonumber\\
&-&\frac{{{\pi ^3}}}{9}\sin \left( {4\pi {x_3}} \right)\cos \left( {2\pi {x_1}} \right) - \frac{{{\pi ^3}}}{9}\sin \left( {2\pi {x_3}} \right)\cos \left( {4\pi {x_1}} \right) +\nonumber\\
&+&\frac{{2{\pi ^3}}}{9}\cos \left( {2\pi {x_3}} \right)\cos \left( {3\pi {x_1}} \right)\sin \left( {3\pi {x_2}} \right) \qquad
\end{eqnarray}\\
%
\subsection{\protect{\bfseries{The coefficients $\chi  _{ij}^{(2)}\left( x \right)$, $i=1,2,3, j=1,11$ and $\mathcal{T} _j\left( t \right)$, $j=3,13$ in $\left(\ref{87}\right)$}}}
%
\begin{eqnarray}
\label{90}
\chi  _{11}^{(2)}\left( x \right) = \sin \left( {\pi {x_1}} \right)\sin \left( {\pi {x_3}} \right) + \cos \left( {\pi {x_1}} \right)\cos \left( {\pi {x_2}} \right)
\end{eqnarray}
\begin{eqnarray}
\label{90}
\chi  _{12}^{(2)}\left( x \right)& = &\cos \left( {\pi {x_1}} \right)\cos \left( {\pi {x_2}} \right)\cos \left( {2\pi {x_3}} \right) - \sin \left( {\pi {x_1}} \right)\cos \left( {2\pi {x_2}} \right)\sin \left( {\pi {x_3}} \right) - \nonumber\\
&-&5\sin \left( {2\pi {x_1}} \right)\sin \left( {\pi {x_2}} \right)\cos \left( {\pi {x_3}} \right)
\end{eqnarray}
\begin{eqnarray}
\label{90}
\chi  _{13}^{(2)}\left( x \right) =  - \sin \left( {2\pi {x_2}} \right)\sin \left( {2\pi {x_3}} \right)
\end{eqnarray}
\begin{eqnarray}
\label{90}
\chi  _{14}^{(2)}\left( x \right)& = &\sin \left( {3\pi {x_1}} \right)\sin \left( {\pi {x_3}} \right) - \cos \left( {3\pi {x_1}} \right)\cos \left( {\pi {x_2}} \right) + \nonumber\\
&+& 3\cos \left( {\pi {x_1}} \right)\cos \left( {3\pi {x_2}} \right) - 3\sin \left( {\pi {x_1}} \right)\sin \left( {3\pi {x_3}} \right)
\end{eqnarray}
\begin{eqnarray}
\label{90}
\chi  _{15}^{(2)}\left( x \right) &= &\cos \left( {3\pi {x_1}} \right)\cos \left( {\pi {x_2}} \right)\cos \left( {2\pi {x_3}} \right) + \cos \left( {\pi {x_1}} \right)\cos \left( {3\pi {x_2}} \right)\cos \left( {2\pi {x_3}} \right) + \nonumber\\
&+& \sin \left( {\pi {x_1}} \right)\cos \left( {2\pi {x_2}} \right)\sin \left( {3\pi {x_3}} \right) - \sin \left( {2\pi {x_1}} \right)\sin \left( {\pi {x_2}} \right)\cos \left( {3\pi {x_3}} \right) + \nonumber\\
&+ &\sin \left( {2\pi {x_1}} \right)\sin \left( {3\pi {x_2}} \right)\cos \left( {\pi {x_3}} \right) + \sin (3\pi {x_1})\cos \left( {2\pi {x_2}} \right)\sin \left( {\pi {x_3}} \right)
\end{eqnarray}
\begin{eqnarray}
\label{90}
\chi  _{16}^{(2)}\left( x \right) =  - \sin \left( {2\pi {x_1}} \right)\cos \left( {2\pi {x_2}} \right) - \sin \left( {2\pi {x_1}} \right)\cos \left( {2\pi {x_3}} \right)
\end{eqnarray}
\begin{eqnarray}
\label{90}
\chi  _{17}^{(2)}\left( x \right) = \cos \left( {3\pi {x_2}} \right)\sin \left( {\pi {x_3}} \right) + \cos \left( {\pi {x_2}} \right)\sin \left( {3\pi {x_3}} \right)
\end{eqnarray}
\begin{eqnarray}
\label{90}
\chi  _{18}^{(2)}\left( x \right) &=& 3\cos \left( {3\pi {x_1}} \right)\sin \left( {2\pi {x_2}} \right)\cos \left( {\pi {x_3}} \right) + 3\sin \left( {3\pi {x_1}} \right)\sin \left( {\pi {x_2}} \right)\sin \left( {2\pi {x_3}} \right) + \nonumber\\
& + &\cos \left( {2\pi {x_1}} \right)\cos \left( {3\pi {x_2}} \right)\sin \left( {\pi {x_3}} \right) - \cos \left( {\pi {x_1}} \right)\sin \left( {2\pi {x_2}} \right)\cos \left( {3\pi {x_3}} \right) - \nonumber\\
&- &\cos \left( {2\pi {x_1}} \right)\cos \left( {\pi {x_2}} \right)\sin \left( {3\pi {x_3}} \right) - \sin \left( {\pi {x_1}} \right)\sin \left( {3\pi {x_2}} \right)\sin \left( {2\pi {x_3}} \right)
\end{eqnarray}
\begin{eqnarray}
\label{90}
\chi  _{19}^{(2)}\left( x \right) =  - \sin \left( {4\pi {x_1}} \right)
\end{eqnarray}
%
\begin{eqnarray}
\label{90}
\chi  _{110}^{(2)}\left( x \right) &=& \sin \left( {2\pi {x_1}} \right)\cos \left( {4\pi {x_3}} \right) + \sin \left( {4\pi {x_1}} \right)\cos \left( {2\pi {x_3}} \right) - \nonumber\\
&-& \sin \left( {4\pi {x_1}} \right)\cos \left( {2\pi {x_2}} \right) - \sin \left( {2\pi {x_1}} \right)\cos \left( {4\pi {x_2}} \right)
\end{eqnarray}
\begin{eqnarray}
\label{90}
\chi  _{111}^{(2)}\left( x \right) = \cos \left( {2\pi {x_1}} \right)\cos \left( {3\pi {x_2}} \right)\sin \left( {3\pi {x_3}} \right)
\end{eqnarray}\\
%
\begin{eqnarray}
\label{90}
\chi  _{21}^{(2)}\left( x \right) = \sin \left( {\pi {x_2}} \right)\sin \left( {\pi {x_1}} \right) + \cos \left( {\pi {x_2}} \right)\cos \left( {\pi {x_3}} \right)
\end{eqnarray}
\begin{eqnarray}
\label{90}
\chi  _{22}^{(2)}\left( x \right)& = &\cos \left( {\pi {x_2}} \right)\cos \left( {\pi {x_3}} \right)\cos \left( {2\pi {x_1}} \right) - \sin \left( {\pi {x_2}} \right)\cos \left( {2\pi {x_3}} \right)\sin \left( {\pi {x_1}} \right) - \nonumber\\
&-&5\sin \left( {2\pi {x_2}} \right)\sin \left( {\pi {x_3}} \right)\cos \left( {\pi {x_1}} \right)
\end{eqnarray}
\begin{eqnarray}
\label{90}
\chi  _{23}^{(2)}\left( x \right) =  - \sin \left( {2\pi {x_3}} \right)\sin \left( {2\pi {x_1}} \right)
\end{eqnarray}
\begin{eqnarray}
\label{90}
\chi  _{24}^{(2)}\left( x \right)& = &\sin \left( {3\pi {x_2}} \right)\sin \left( {\pi {x_1}} \right) - \cos \left( {3\pi {x_2}} \right)\cos \left( {\pi {x_3}} \right) + \nonumber\\
&+& 3\cos \left( {\pi {x_2}} \right)\cos \left( {3\pi {x_3}} \right) - 3\sin \left( {\pi {x_2}} \right)\sin \left( {3\pi {x_1}} \right)
\end{eqnarray}
\begin{eqnarray}
\label{90}
\chi  _{25}^{(2)}\left( x \right) &= &\cos \left( {3\pi {x_2}} \right)\cos \left( {\pi {x_3}} \right)\cos \left( {2\pi {x_1}} \right) + \cos \left( {\pi {x_2}} \right)\cos \left( {3\pi {x_3}} \right)\cos \left( {2\pi {x_1}} \right) + \nonumber\\
&+& \sin \left( {\pi {x_2}} \right)\cos \left( {2\pi {x_3}} \right)\sin \left( {3\pi {x_1}} \right) - \sin \left( {2\pi {x_2}} \right)\sin \left( {\pi {x_3}} \right)\cos \left( {3\pi {x_1}} \right) + \nonumber\\
&+ &\sin \left( {2\pi {x_2}} \right)\sin \left( {3\pi {x_3}} \right)\cos \left( {\pi {x_1}} \right) + \sin (3\pi {x_2})\cos \left( {2\pi {x_3}} \right)\sin \left( {\pi {x_1}} \right)
\end{eqnarray}
\begin{eqnarray}
\label{90}
\chi  _{26}^{(2)}\left( x \right) =  - \sin \left( {2\pi {x_2}} \right)\cos \left( {2\pi {x_3}} \right) - \sin \left( {2\pi {x_2}} \right)\cos \left( {2\pi {x_1}} \right)
\end{eqnarray}
\begin{eqnarray}
\label{90}
\chi  _{27}^{(2)}\left( x \right) = \cos \left( {3\pi {x_3}} \right)\sin \left( {\pi {x_1}} \right) + \cos \left( {\pi {x_3}} \right)\sin \left( {3\pi {x_1}} \right)
\end{eqnarray}
\begin{eqnarray}
\label{90}
\chi  _{28}^{(2)}\left( x \right) &=& 3\cos \left( {3\pi {x_2}} \right)\sin \left( {2\pi {x_3}} \right)\cos \left( {\pi {x_1}} \right) + 3\sin \left( {3\pi {x_2}} \right)\sin \left( {\pi {x_3}} \right)\sin \left( {2\pi {x_1}} \right) + \nonumber\\
& + &\cos \left( {2\pi {x_2}} \right)\cos \left( {3\pi {x_3}} \right)\sin \left( {\pi {x_1}} \right) - \cos \left( {\pi {x_2}} \right)\sin \left( {2\pi {x_3}} \right)\cos \left( {3\pi {x_1}} \right) - \nonumber\\
&- &\cos \left( {2\pi {x_2}} \right)\cos \left( {\pi {x_3}} \right)\sin \left( {3\pi {x_1}} \right) - \sin \left( {\pi {x_2}} \right)\sin \left( {3\pi {x_3}} \right)\sin \left( {2\pi {x_1}} \right)
\end{eqnarray}
\begin{eqnarray}
\label{90}
\chi  _{29}^{(2)}\left( x \right) =  - \sin \left( {4\pi {x_2}} \right)
\end{eqnarray}
%
\begin{eqnarray}
\label{90}
\chi  _{210}^{(2)}\left( x \right) &=& \sin \left( {2\pi {x_2}} \right)\cos \left( {4\pi {x_1}} \right) + \sin \left( {4\pi {x_2}} \right)\cos \left( {2\pi {x_1}} \right) - \nonumber\\
&-& \sin \left( {4\pi {x_2}} \right)\cos \left( {2\pi {x_3}} \right) - \sin \left( {2\pi {x_2}} \right)\cos \left( {4\pi {x_3}} \right)
\end{eqnarray}
\begin{eqnarray}
\label{90}
\chi  _{211}^{(2)}\left( x \right) = \cos \left( {2\pi {x_2}} \right)\cos \left( {3\pi {x_3}} \right)\sin \left( {3\pi {x_1}} \right)
\end{eqnarray}\\
%
\begin{eqnarray}
\label{90}
\chi  _{31}^{(2)}\left( x \right) = \sin \left( {\pi {x_3}} \right)\sin \left( {\pi {x_2}} \right) + \cos \left( {\pi {x_3}} \right)\cos \left( {\pi {x_1}} \right)
\end{eqnarray}
\begin{eqnarray}
\label{90}
\chi  _{32}^{(2)}\left( x \right)& = &\cos \left( {\pi {x_3}} \right)\cos \left( {\pi {x_1}} \right)\cos \left( {2\pi {x_2}} \right) - \sin \left( {\pi {x_3}} \right)\cos \left( {2\pi {x_1}} \right)\sin \left( {\pi {x_2}} \right) - \nonumber\\
&-&5\sin \left( {2\pi {x_3}} \right)\sin \left( {\pi {x_1}} \right)\cos \left( {\pi {x_2}} \right)
\end{eqnarray}
\begin{eqnarray}
\label{90}
\chi  _{33}^{(2)}\left( x \right) =  - \sin \left( {2\pi {x_1}} \right)\sin \left( {2\pi {x_2}} \right)
\end{eqnarray}
\begin{eqnarray}
\label{90}
\chi  _{34}^{(2)}\left( x \right)& = &\sin \left( {3\pi {x_3}} \right)\sin \left( {\pi {x_2}} \right) - \cos \left( {3\pi {x_3}} \right)\cos \left( {\pi {x_1}} \right) + \nonumber\\
&+& 3\cos \left( {\pi {x_3}} \right)\cos \left( {3\pi {x_1}} \right) - 3\sin \left( {\pi {x_3}} \right)\sin \left( {3\pi {x_2}} \right)
\end{eqnarray}
\begin{eqnarray}
\label{90}
\chi  _{35}^{(2)}\left( x \right) &= &\cos \left( {3\pi {x_3}} \right)\cos \left( {\pi {x_1}} \right)\cos \left( {2\pi {x_2}} \right) + \cos \left( {\pi {x_3}} \right)\cos \left( {3\pi {x_1}} \right)\cos \left( {2\pi {x_2}} \right) + \nonumber\\
&+& \sin \left( {\pi {x_3}} \right)\cos \left( {2\pi {x_1}} \right)\sin \left( {3\pi {x_2}} \right) - \sin \left( {2\pi {x_3}} \right)\sin \left( {\pi {x_1}} \right)\cos \left( {3\pi {x_2}} \right) + \nonumber\\
&+ &\sin \left( {2\pi {x_3}} \right)\sin \left( {3\pi {x_1}} \right)\cos \left( {\pi {x_2}} \right) + \sin (3\pi {x_3})\cos \left( {2\pi {x_1}} \right)\sin \left( {\pi {x_2}} \right)
\end{eqnarray}
\begin{eqnarray}
\label{90}
\chi  _{36}^{(2)}\left( x \right) =  - \sin \left( {2\pi {x_3}} \right)\cos \left( {2\pi {x_1}} \right) - \sin \left( {2\pi {x_3}} \right)\cos \left( {2\pi {x_2}} \right)
\end{eqnarray}
\begin{eqnarray}
\label{90}
\chi  _{37}^{(2)}\left( x \right) = \cos \left( {3\pi {x_1}} \right)\sin \left( {\pi {x_2}} \right) + \cos \left( {\pi {x_1}} \right)\sin \left( {3\pi {x_2}} \right)
\end{eqnarray}
\begin{eqnarray}
\label{90}
\chi  _{38}^{(2)}\left( x \right) &=& 3\cos \left( {3\pi {x_3}} \right)\sin \left( {2\pi {x_1}} \right)\cos \left( {\pi {x_2}} \right) + 3\sin \left( {3\pi {x_3}} \right)\sin \left( {\pi {x_1}} \right)\sin \left( {2\pi {x_2}} \right) + \nonumber\\
& + &\cos \left( {2\pi {x_3}} \right)\cos \left( {3\pi {x_1}} \right)\sin \left( {\pi {x_2}} \right) - \cos \left( {\pi {x_3}} \right)\sin \left( {2\pi {x_1}} \right)\cos \left( {3\pi {x_2}} \right) - \nonumber\\
&- &\cos \left( {2\pi {x_3}} \right)\cos \left( {\pi {x_1}} \right)\sin \left( {3\pi {x_2}} \right) - \sin \left( {\pi {x_3}} \right)\sin \left( {3\pi {x_1}} \right)\sin \left( {2\pi {x_2}} \right)
\end{eqnarray}
\begin{eqnarray}
\label{90}
\chi  _{39}^{(2)}\left( x \right) =  - \sin \left( {4\pi {x_3}} \right)
\end{eqnarray}
%
\begin{eqnarray}
\label{90}
\chi  _{310}^{(2)}\left( x \right) &=& \sin \left( {2\pi {x_3}} \right)\cos \left( {4\pi {x_2}} \right) + \sin \left( {4\pi {x_3}} \right)\cos \left( {2\pi {x_2}} \right) - \nonumber\\
&-& \sin \left( {4\pi {x_3}} \right)\cos \left( {2\pi {x_1}} \right) - \sin \left( {2\pi {x_3}} \right)\cos \left( {4\pi {x_1}} \right)
\end{eqnarray}
\begin{eqnarray}
\label{90}
\chi  _{311}^{(2)}\left( x \right) = \cos \left( {2\pi {x_3}} \right)\cos \left( {3\pi {x_1}} \right)\sin \left( {3\pi {x_2}} \right)
\end{eqnarray}\\
\begin{eqnarray}
\label{90}
{\mathcal{T} _3}\left( t \right) = \frac{{{e^{ - 2{\pi ^2}\kappa t}} - 3{e^{ - 6{\pi ^2}\kappa t}} + 2{e^{ - 8{\pi ^2}\kappa t}}}}{{96{\left(\pi \kappa\right) ^2}}}
\end{eqnarray}
\begin{eqnarray}
\label{90}
{\mathcal{T} _4}\left( t \right) = \frac{{2{\pi ^2}\left( {\kappa t} \right){e^{ - 6{\pi ^2}\kappa t}} + {e^{ - 6{\pi ^2}\kappa t}} - {e^{ - 8{\pi ^2}\kappa t}}}}{{12{\left(\pi \kappa\right) ^2}}},\nonumber\\  \protect{\textrm{In the first term in the numerator}}\,\, \varepsilon=6,\,\, \sigma=1
\end{eqnarray}
%
\begin{eqnarray}
\label{90}
{\mathcal{T} _5}\left( t \right) = \frac{{{e^{ - 6{\pi ^2}\kappa t}} - {e^{ - 8{\pi ^2}\kappa t}} - 2{\pi ^2}\left( {\kappa t} \right){e^{ - 8{\pi ^2}\kappa t}}}}{{6{\left(\pi \kappa\right) ^2}}},\nonumber\\  \protect{\textrm{In the last term in the numerator}}\,\, \varepsilon=8,\,\, \sigma=1
\end{eqnarray}
%
\begin{eqnarray}
\label{90}
{\mathcal{T} _6}\left( t \right) = \frac{{{e^{ - 6{\pi ^2}\kappa t}} - 2{e^{ - 8{\pi ^2}\kappa t}} + {e^{ - 10{\pi ^2}\kappa t}}}}{{48{\left(\pi \kappa\right) ^2}}}
\end{eqnarray}
%
\begin{eqnarray}
\label{90}
{\mathcal{T} _7}\left( t \right) = \frac{{3{e^{ - 6{\pi ^2}\kappa t}} - 4{e^{ - 8{\pi ^2}\kappa t}} + {e^{ - 14{\pi ^2}\kappa t}}}}{{288{\left(\pi \kappa\right) ^2}}}
\end{eqnarray}
%
\begin{eqnarray}
\label{90}
{\mathcal{T} _8}\left( t \right) = \frac{{4{\pi ^2}\left( {\kappa t} \right){e^{ - 8{\pi ^2}\kappa t}} - 3{e^{ - 8{\pi ^2}\kappa t}} + 4{e^{ - 10{\pi ^2}\kappa t}} - {e^{ - 12{\pi ^2}\kappa t}}}}{{72{\left(\pi \kappa\right) ^3}}},\nonumber\\  \protect{\textrm{In the first term in the numerator}}\,\,\varepsilon=8,\,\, \sigma=1
\end{eqnarray}
%
\begin{eqnarray}
\label{90}
{\mathcal{T} _9}\left( t \right) = \frac{{{e^{ - 8{\pi ^2}\kappa t}} - 4{\pi ^2}\left(\kappa t\right){e^{ - 10{\pi ^2}\kappa t}} - {e^{ - 12{\pi ^2}\kappa t}}}}{{24{\left(\pi \kappa\right) ^3}}},\nonumber\\  \protect{\textrm{In the second term in the numerator}}\,\, \varepsilon=10,\,\, \sigma=1
\end{eqnarray}
%
\begin{eqnarray}
\label{90}
{\mathcal{T} _{10}}\left( t \right) = \frac{{{e^{ - 8{\pi ^2}\kappa t}} - 3{e^{ - 10{\pi ^2}\kappa t}} + 3{e^{ - 12{\pi ^2}\kappa t}} - {e^{ - 14{\pi ^2}\kappa t}}}}{{216{\left(\pi \kappa\right) ^3}}}
\end{eqnarray}
%
\begin{eqnarray}
\label{90}
{\mathcal{T} _{11}}\left( t \right) = \frac{{3{e^{ - 8{\pi ^2}\kappa t}} - 8{e^{ - 10{\pi ^2}\kappa t}} + 6{e^{ - 12{\pi ^2}\kappa t}} - {e^{ - 16{\pi ^2}\kappa t}}}}{{432{\left(\pi \kappa\right) ^3}}}
\end{eqnarray}
%
\begin{eqnarray}
\label{90}
{\mathcal{T} _{12}}\left( t \right) = \frac{{10{e^{ - 8{\pi ^2}\kappa t}} - 24{e^{ - 10{\pi ^2}\kappa t}} - {e^{ - 20{\pi ^2}\kappa t}} + 15{e^{ - 12{\pi ^2}\kappa t}}}}{{4320{\left(\pi \kappa\right) ^3}}}
\end{eqnarray}
%
\begin{eqnarray}
\label{90}
{\mathcal{T} _{13}}\left( t \right) = \frac{{15{e^{ - 8{\pi ^2}\kappa t}} - 35{e^{ - 10{\pi ^2}\kappa t}} - {e^{ - 22{\pi ^2}\kappa t}} + 21{e^{ - 12{\pi ^2}\kappa t}}}}{{3780{\left(\pi \kappa\right) ^3}}}
\end{eqnarray}
\end{document}